\documentclass[a4paper,twoside, draft, 11pt]{amsart}

\usepackage{amssymb,amsfonts,amsmath,amsthm,amscd,verbatim,epic,pstricks}

\newtheorem{theor}{~~~~Theorem}
\newtheorem{prop}{~~~~Proposition}
\newtheorem{cor}{~~~~Corollary}
\newtheorem{lemma}{~~~~Lemma}

\newtheorem{stat}{~~~~ Statement}
\newtheorem{remark}{~~~~Remark}
\newtheorem{defin}{~~~~Definition}

\newcommand\halmos{\rule{0.1in}{0.1in}}

\begin{document}
\title{
Differential geometry of curves in Lagrange Grassmannians with given
Young diagram}
\date{}
\author{
Igor Zelenko\address{S.I.S.S.A., Via Beirut 2-4, 34014, Trieste,
Italy; email: zelenko@sissa.it} \and Chengbo Li \address{S.I.S.S.A.,
Via Beirut 2-4, 34014, Trieste, Italy; email: chengbo@sissa.it}}
\maketitle \markboth {Igor Zelenko and Chengbo Li} {Differential
geometry of curves in Lagrange Grassmannians}
\begin{abstract}
Curves in Lagrange Grassmannians appear naturally in the intrinsic study of
geometric structures on manifolds. By a smooth geometric structure
on a manifold  we mean any submanifold of its tangent bundle,
transversal to the fibers. One can consider the time-optimal problem
naturally associate with a geometric structure. The Pontryagin
extremals of this optimal problem are integral curves of certain
Hamiltonian system in the cotangent bundle. The dynamics of the
fibers of the cotangent bundle w.r.t. this system along an extremal
is described by certain curve in a Lagrange Grassmannian, called
Jacobi curve of the extremal. Any symplectic invariant of the Jacobi
curves produces the invariant of the original geometric structure.
The basic characteristic of a curve in a Lagrange Grassmannian is
its Young diagram. The number of boxes in its $k$th column is equal
to the rank of the $k$th derivative of the curve (which is an
appropriately defined linear mapping) at a generic point.
We will describe the construction  of the complete system of
symplectic invariants for parameterized curves in a Lagrange
Grassmannian with given Young diagram. It allows to develop in a
unified way local differential geometry of very wide classes of
geometric structures on manifolds, including both classical
geometric structures such as Riemannian and Finslerian structures
and less classical ones such as sub-Riemannian and sub-Finslerian
structures, defined on nonholonomic distributions.
\end{abstract}
\section{Introduction}

Let $W$ be a $2m$-dimensional linear space endowed  with a
symplectic form $\omega$. Recall that an $m$-dimensional subspace
$\Lambda$ of $W$ is called Lagrangian, if $\omega|_\Lambda=0$.
Lagrange Grassmannian $L(W)$ of $W$
 is the  set of all Lagrangian subspaces of $W$.
The linear Symplectic group acts naturally on $L(W)$.
Invariants of curves in a Lagrange Grassmannian w.r.t. this action
are called symplectic. The present paper is devoted to the
construction of a complete system of symplectic invariants for
smooth parameterized curves in the Lagrange Grassmannian $L(W)$, i.e., a
set of invariants (independent one of each other) such that there
exists the unique, up to a symplectic transformation, curve in
$L(W)$ with the prescribed invariants from this set. Of course, this
problem is a particular case of the classical problem on
differential geometry of curves in homogeneous spaces. The general
procedure for the latter problem was developed already by E. Cartan
with his method of moving frames. On the other hand,
by studying
curves in Lagrange Grasmannians, one can develop in a unified way
local differential geometry of very wide classes of geometric structures
on manifolds, including both classical geometric structures such as
Riemannian and Finslerian structures and less classical such as
sub-Riemannian or sub-Finslerian structures.\footnote{Differential
geometry of rank 2 vector distributions (without additional
structures on them) can be treated as well by studying
unparameterized curves in Lagrange Grassmannians
(\cite{doubzel},\cite{doubzel1})}. Therefore, the explicit
construction of moving frames and invariants in the particular
situation of curves in Lagrange Grassmannians
is important by itself.


Let us briefly  describe how curves in Lagrange Grassmannians appear
in intrinsic study of geometric structures (more detailed and
general presentation can be found in
\cite{agrgam1} or \cite{jac1}). Here by a smooth geometric structure
on a manifold $M$ we mean any submanifold $\mathcal V\subset TM$,
transversal to fibers.
Let $\mathcal V_q=\mathcal V\cap T_qM$. For example,
if $\mathcal V_q$ is an intersection of an ellipsoid centered at the
origin with a linear subspace ${\mathcal D}_q$ in $T_qM$ (where both
the ellipsoids and the subspaces $\mathcal D_q$ depend smoothly on $q$),
then $\mathcal V$ is called a \emph {sub-Riemannian structure on
$M$ with underlying distribution $\mathcal D$}. In this
case $\mathcal V_q$ is the unit sphere w.r.t. the unique Euclidean
norm $||\cdot||_q$ on $\mathcal D_q$, i.e. fixing an ellipsoid in
$\mathcal D_q$ is equivalent to fixing an Euclidean norm on
$\mathcal D_q$ for any $q\in M$. This reformulation justifies the
term ``sub-Riemannian''. In particular, it defines in the obvious way
the length of any curve tangent to the underlying distribution.  If
in the constructions above we replace the ellipsoids by the
boundaries of strongly convex bodies in $T_qM$
(sometimes also assumed to be  symmetric w.r.t. the origin)
we will get a sub-Finslerian structure on $M$. Note also
that, if the underlying distribution $\mathcal D=TM$, we get just a
Riemannian (a Finslerian) structure on $M$.

Actually, one can  look at a geometric structure $\mathcal V$ as a
control system on $M$: the set $\mathcal V_q$ defines the set of all
admissible velocities of motion from the point $q$. A Lipshitzian
curve $\gamma:[0,T]\mapsto M$ is called an \emph {admissible
trajectory of $\mathcal V$}, if $\dot \gamma(t)\in \mathcal
V_\gamma(t)$ for a. e. $t$. Now one can consider the time-optimal
problem on $\mathcal V$: given two points $q_0$ and $q_1$ to find an
admissible trajectory, steering from $q_0$ to $q_1$ in a minimal
time. The extremals of this optimal problem are obtained from the
Pontryagin Maximum Principle of Optimal Control Theory
(\cite{pontr}). Here for simplicity of presentation let us suppose
that the maximized Hamiltonian of the Pontryagin Maximum Principle
\begin{equation}
\label{maxH} H(p, q)=\max_{v\in \mathcal V_q} p(v), \quad q\in M,
p\in T_q^*M
\end{equation}
is well defined and smooth in an open domain $O\subset T^*M$ and for
some $c>0$ (and therefore for any $c>0$ by homogeneity of $H$ on
each fiber of $T^*M$) the corresponding level set $$\mathcal
H_c=\{\lambda\in O: H(\lambda)=c\}$$ is nonempty and consists of
regular points of $H$. Consider the Hamiltonian vector field $\vec
H$ on $\mathcal H_c$, corresponding to the Hamiltonian $H$, i.e. the
vector field satisfying $i_{\vec H}\bar \omega=-dH$, where
$\bar\omega$ is the canonical symplectic structure on $T^*M$. The
integral curves of this Hamiltonian system are normal Pontryagin
extremals of the time-optimal problem, associated with geometric
structure $\mathcal V$, or, shortly, normal extremals of $\mathcal
V$. For example, if $\mathcal V$ is a sub-Riemannian structure with
underlying distribution $\mathcal D$, then the maximized Hamiltonian
satisfies $H(p,q)=||p|_{_{\mathcal D_q}}||_q$, i.e. $H(p,q)$ is
equal to the norm of the restriction of the functional $p\in T_q^*M$
on $\mathcal D_q$ w.r.t. the Euclidean norm $||\cdot||_q$ on
$\mathcal D_q$; $O=T^*M\backslash \mathcal D^\perp$, where $\mathcal
D^\perp$ is the annihilator of $\mathcal D$,
$$\mathcal D^\perp=\{(p,q)\in T^*M: p(v)=0\,\,\forall v\in \mathcal D_q\}.$$ The
projections of the trajectories of the corresponding Hamiltonian
systems to the base manifold $M$ are normal sub-Riemannian
geodesics. If $\mathcal D=TM$, then they are exactly the Riemannian
geodesics of the corresponding Riemannian structure.

Further let $\mathcal H_c(q)=\mathcal H_c\cap T_q^*M$. $\mathcal
H_c(q)$ is a codimension 1 submanifold of $T^*_qM$.  For any
$\lambda\in \mathcal H_c$ denote $\Pi_\lambda=T_\lambda
\bigr(\mathcal H_c(\pi(\lambda))\bigr)$, where $\pi:T^*M\mapsto M$
is the canonical projection. Actually $\Pi_\lambda$ is the vertical
subspace of $T_\lambda \mathcal H_c$,
\begin{equation}
\label{pilam} \Pi_\lambda=\{\xi\in T_\lambda \mathcal
H_c:\pi_*(\xi)=0\}.
\end{equation}

Now with any integral curve of $\vec H$ one can associate a curve in
a Lagrange Grassmannian, which describes the dynamics of the
vertical subspaces $\Pi_\lambda$ along this integral curve  w.r.t.
the flow $e^{t\vec H}$, generated by $\vec H$.  For this let
\begin{equation}
t\mapsto J_\lambda(t)\stackrel{def}{=}
e^{-t \vec{H}}_* \bigl(\Pi_{e^{t \vec H}\lambda}\bigr)
\Bigr/\{\mathbb R\vec H(\lambda)\}\Bigl..\label{jcurve}
\end{equation}
The curve  $J_\lambda(t)$ is the curve in the Lagrange Grassmannian
of the linear symplectic space $W_\lambda=T_\lambda \mathcal
H_c/\{\mathbb R\vec H(\lambda)\}$ (endowed with the symplectic form
$\omega$ induced in the obvious way by the canonical symplectic form
$\bar\omega$ of $T^*M$). It is called {\sl the Jacobi curve of the
curve $ e^{t\vec{\mathcal H}}\lambda$
 attached at the point $\lambda$}. Note  also that if
$\bar\lambda=e^{\bar t\vec H}\lambda$  and $\Phi: W_\lambda\mapsto
W_{\bar\lambda}$ is a symplectic transformation induced in the
natural way by a linear mapping $e^{t\vec H}_*:T_\lambda \mathcal
H_c\mapsto T_{\bar\lambda}\mathcal H_c$, then by (\ref{jcurve}) we
have
\begin{equation}
\label{diffpoint} J_{\bar\lambda}(t)=\Phi\bigl(J_ \lambda(t-\bar
t)\bigr).
\end{equation}
In other words, the Jacobi curves of the same integral curve of
$\vec H$ attached at two different points of this curve are the
same, up to symplectic transformation between the corresponding
ambient linear symplectic spaces and the corresponding shift of the
parameterizations. Therefore, any symplectic invariant of the Jacobi
curve produces the function  on the manifold $\mathcal H_c$,
intrinsically related to the geometric structure $V$ (the value of
this function at $\lambda\in\mathcal H_c$ is equal to the value of
the chosen symplectic invariant of the curve $J_\lambda(t)$ at
$t=0$). In this way the problem of finding differential invariants
of geometric structure can be  essentially reduced to the much more
treatable problem of finding symplectic invariants of certain curves
in a Lagrange Grassmannian.

In all constructions above one can replace the
maximized Hamiltonian $H$ by some its power $H^s$. It causes only the
reparametrization of the Jacobi curve of the type $t\mapsto Ct$ for
some constant $C$. For example in the case of sub-Riemannian
structures it is more convenient to work with $H^2$ instead of $H$,
because $H^2$ is a polynomial on the fibers of $T^*M$.

Jacobi curves of integral curves of $\vec H$ are not arbitrary
curves of Lagrangian Grassmannian but they inherit special features
of the geometric structure $\mathcal V$. To specify these features
recall that the tangent space $T_\Lambda L(W)$ to the Lagrangian
Grassmannian $L(W)$ at the point $\Lambda$ can be naturally
identified with the space ${\rm Quad}(\Lambda)$ of all quadratic
forms on linear space $\Lambda\subset W$.
Namely, given $\mathfrak V\in T_\Lambda L(W)$ take a curve
$\Lambda(t)\in L(W)$ with $\Lambda(0)=\Lambda$ and
$\dot\Lambda=\mathfrak V$. Given some vector $l\in\Lambda$, take a
curve $\ell(\cdot)$ in $W$ such that $\ell(t)\in \Lambda(t)$ for all
$t$ and $\ell(0)=l$. Define the quadratic form
\begin{equation}
\label{quad} Q_{\mathfrak V}(l)=\omega(\frac{d}{dt}\ell(0),l).
\end{equation}
Using the fact that the spaces $\Lambda(t)$ are Lagrangian,
it is easy to see that $Q_{\mathfrak V}(l)$ does not depend on the
choice of the curves $\ell$  and $\Lambda(t)$ with the above
properties, but depends only on $\mathfrak V$.
So, we have the linear mapping from $T_\Lambda L(W)$ to the spaces
${\rm Quad}(\Lambda)$, $\mathfrak V\mapsto Q_{\mathfrak V}$.
A simple counting of dimensions shows that this mapping is a
bijection and it defines the required identification. A curve
$\Lambda(\cdot)$ in a Lagrange Grassmannian is called \emph{regular
at a point $\tau$}, if its velocity at $\tau$ is a nondegenerated
quadratic form, and \emph{nonregular at $\tau$} otherwise. The rank
of the velocity $ \dot\Lambda(\tau)$ of a curve $\Lambda(\cdot)$ at
a point $\tau$ is  called shortly the \emph {rank of
$\Lambda(\cdot)$} at $\tau$. A curve $\Lambda(\cdot)$ is called
\emph{monotonically nondecreasing (nonincreasing)} if the velocity
is nonnegative (nonpositive) definite at any point. We also will
call such curves \emph{monotonic}.

It turns out (see, for example, \cite[Proposition 1]{jac1}) that the
velocity of the Jacobi curve $J_\lambda(t)$ at $t=0$ is
equal to the restriction of the Hessian of $H$ to the tangent space
to $\mathcal H_{H(\lambda)}$ at $\lambda$. This together with
\eqref{diffpoint} implies easily (\cite{jac1}) that the rank of the
Jacobi curve $J_\lambda(t)$ at $t=\tau$ is not greater then $\dim\,
\mathcal V_{\pi(e^{\tau\vec H}\lambda)}$. For sub-Riemannian
structures the rank of Jacobi curves at any point is equal to $\text
{rank}\, \mathcal D-1$, where $\mathcal D$ is the underlying
distribution, i.e., except the case $\mathcal D=TM$ (corresponding
to a Riemannian structure), the Jacobi curves appearing in
sub-Riemannian structures are nonregular at any point.

Regular curves were treated in \cite {agrgam1}, where the notion of
the curvature operator was introduced (the work \cite{ovs} is
closely related as well). In particular, calculating the curvature
operator for Jacobi curves, associated with a Riemannian structure,
one gets a part of the Riemannian curvature tensor, appearing in the
classical Jacobi equation for Jacobi vector fields along the
Riemannian geodesics. The whole Riemannian curvature tensor  can be
recovered uniquely from it.

Basic symplectic invariants of curves (both parameterized and
unparameterized) in Lagrange Grassmannians, which are nonregular at
any point,  were constructed in \cite{jac1}, using the notion of
cross-ratio of four points in Lagrange Grassmannians. But the only
nonregular (at any point) curves in Lagrange Grassmannians, for
which the complete system of symplectic invariants was constructed,
were parameterized curves of constant rank 1 (\cite{zel}).

In the present paper we develop differential geometry of curves of
any constant rank in Lagrange Grassmannians, implementing the scheme
briefly described in the Introduction of \cite{zel}. In the study of
generic germs of nonregular curves the basic characteristic are not
only the rank of its velocity, but a certain Young diagram (see
subsections 2.1). The rank of the curve is the number of boxes in
the first column of this Young diagram. It is also very convenient to
consider the additional "smaller ' diagram, called the reduced Young
diagram (see subsections 2.2).
For a regular curve the Young diagram consists of one column and the
reduced Young diagram consists of one box, while for rank 1 curve
the Young diagram and its reduction coincide and consist of one row.
For any monotonic curve or a generic nonmonotonic curve
$\Lambda(\cdot)$ in a Lagrange Grassmannian with given Young diagram
we construct the principal bundle (over the curve itself) of frames
in the ambient symplectic space endowed with the canonical principal
connection or the bundle of moving frames, canonically associated
with the curve
(Theorems \ref{maintheor} and \ref{maintheor1}). These moving frames
are defined by the form of the matrix in their structural equation.
During the process of normalization we get the canonical splitting
of any subspaces $\Lambda(t)$ such that the subspaces of the
splitting are parameterized by boxes of the reduced Young diagram
and each subspace of the splitting is endowed with the canonical
Euclidean or pseudo-Euclidean structure (in the monotonic and
nonmonotonic cases repsectively). Also we construct in a canonical
way the additional curve $\Lambda^{trans}(\cdot)$ in a Lagrange
Grassmannian such that any subspace $\Lambda^{trans}(t)$ is
transversal to the subspace $\Lambda(t)$ for any $t$. Further, using
the matrix in the structural equation of canonical moving frames, we
obtain the tuple of one-parametic families of linear mappings
between the subspaces of the canonical splitting. This tuple
constitutes a kind of a complete system of symplectic invariants of
the curve in a sense formulated  in terms of quivers and their
representations (Theorems \ref{size1} and \ref{size2}).
In the case when the Young diagram of the curve $\Lambda(\cdot)$
has no rows with the same number of boxes, we get in this way a
complete system of scalar invariants of the curve $\Lambda(\cdot)$
in the usual sense.
As the consequences of our constructions  in section 5 we get the
canonical (non-linear) connection on an open sets of the cotangent
bundle $T^*M$, the curvature-type invariants, and additional
nontrivial structures on the fibers of $T^*M$ for geometric
structures on a manifold $M$, including sub-Riemannian and
sub-Finslerian structures, satisfying very general assumptions.

\section{The main results}
\setcounter{equation}{0}
\subsection{
The flag and the Young diagrams associated with a curve}
With any curve $\Lambda(\cdot)$ in Grassmannian $G_k(W)$ of
$k$-dimensional subspaces of a linear space $W$ one can associate a
curve of flags of subspaces in $W$. For this let $\mathfrak
S(\Lambda)$ be the set of all smooth curves  $\ell(t)$ in $W$ such
that $\ell(t)\in \Lambda(t)$ for all $t$.
Denote
\begin{equation}
\label{primeiext} \Lambda^{(i)}(\tau)={\rm span}\Bigl\{
\frac{d^j}{d\tau^j}\ell(\tau)|: \ell\in \mathfrak S(\Lambda),\,
0\leq j\leq i\Bigr\}.
\end{equation}

The subspaces $\Lambda^{(i)}(\tau)$
are  called {\it the $i$th extension}
of the curve $\Lambda(\cdot)$ at the point $\tau$. Recall that the
tangent space $T_\Lambda G_k(W)$ to any subspace $\Lambda\in G_k(W)$
can be identified with the space ${\rm Hom}\, (\Lambda,W/\Lambda)$
of linear mappings from $\Lambda$ to $W/\Lambda$. Using this
identification, if $P:\Lambda\mapsto W/\Lambda$ is the canonical
projection to the factor, then
$\Lambda^{(1)}(\tau)=P^{(-1)}\bigl({\rm
Im}\,\dot\Lambda(\tau)\bigr)$, which implies that $\dim
\Lambda^{(1)}(\tau)- \dim \Lambda(\tau)=\text{rank}\,
\dot\Lambda(\tau)$. By construction
$\Lambda^{(i-1)}(\tau)\subseteq\Lambda^{(i)}(\tau)$.
The flag
\begin{equation}
\label{flag}
\Lambda(\tau)\subseteq\Lambda^{(1)}(\tau)
\subseteq\Lambda^{(2)}(\tau) \subseteq\ldots
\end{equation}
is called the {\it associated (right) flag of the curve
$\Lambda(\cdot)$ at the point $t$}.

From now on we suppose that dimensions of all subspaces
$\Lambda^{(i)}(t)$ (and therefore of $\Lambda_{(i)}(t)$) are
independent of $t$. In this case from (\ref{primeiext}) it is easy
to obtain that the following inequalities hold
\begin{equation}
\label{ineq}
\dim \Lambda^{(i+1)}-\dim \Lambda^{(i)}\leq \dim \Lambda^{(i)}-\dim
\Lambda^{(i-1)}.
\end{equation}
Using inequalities \eqref{ineq}, to any curve $\Lambda(\cdot)$ we
can assign the Young diagram in the following way: the number of
boxes in the $i$th column of this Young diagram is equal to $\dim
\Lambda^{(i)}-\dim \Lambda^{(i-1)}$. It will be called the {\it
Young diagram of the curve $\Lambda(\cdot)$}. In particular, the
number of boxes in the first column is equal to the rank of the
curve.

Now
suppose that $W$ is an even-dimensional linear space endowed with a
symplectic structure $\omega$ and the curve $\Lambda(\cdot)$ is a
curve in the  Lagrangian Grassmannian  $L(W)$.

\begin{remark}
\label{loss} {\rm Without loss of generality, we will suppose that
there exists an integer $p$ such that $\Lambda^{(p)}(t)=W$.
Otherwise, if $\Lambda^{(p+1)}(t)=\Lambda^{(p)}(t)\subsetneq W$, then
the subspace $\Lambda^{(p)}(t)$ does not depend on $t$. Set $V=
\Lambda^{(p)}(t)$. Then $V^\angle\subset \Lambda(t)$ for any $t$ and
all information about the original curve $\Lambda(\cdot)$ is
contained in the curve $\Lambda(\cdot)/V^\angle$, which is the curve
of Lagrangian subspaces in the symplectic space $V/V^\angle$, and
the $p$th extension of the curve $\Lambda(\cdot)/V^\angle$ is equal
to $V/V^\angle$. So, we can work with the curve
$\Lambda(\cdot)/V^\angle$ and the symplectic space $V/V^\angle$
instead of the curve $\Lambda(\cdot)$ and the symplectic space $W$.}
\end{remark}



\subsection{ The normal moving frame.}
The Young diagram is a basic invariant of the curve in Lagrange
Grassmannians. As indices of vectors in our Darboux moving frames we
will take the boxes of the Young diagram instead of the natural
numbers. We found it extremely useful both for formulation of our
results and their proofs.

First note that any Young diagram $D$  can be uniquely represented
as a union of $d$ rectangular diagrams $D_i$ of the sizes $r_i\times
p_i$, $1\leq i\leq d$, such that the sequence $\{p_i\}_{i=1}^d$ is
strictly decreasing. The Young diagram $\Delta$, consisting of $d$
rows such that the $i$th row has $p_i$ boxes, will be called the
{\emph reduced diagram} or the \emph{reduction of the diagram D}.
In order to distinguish between boxes and rows of the diagram $D$
and its reduction $\Delta$, the boxes of $\Delta$ will be called
\emph{superboxes} and the rows of $\Delta$ will be called
\emph{levels}. To the $j$th superbox $a$ of the $i$th level of
$\Delta$ one can assign the $j$th column of the rectangular subdiagram
$D_i$ of $D$ and the integer number $r_i$ (equal to the number of
boxes in this subcolumn), called the \emph {size} of the superbox
$a$.

As usual, by $\Delta\times \Delta$ we will mean the set of pairs of
superboxes of $\Delta$. Also denote by $\rm {Mat}$ the set of
matrices of all sizes. The mapping $R: \Delta\times \Delta\mapsto
\rm {Mat}$
is called \emph {compatible with the Young diagram $D$},
 if to any pair $(a, b)$ of superboxes of sizes $s_1$ and $s_2$ respectively
the matrix $R(a,b)$ is of the size $s_2\times s_1$. The compatible
mapping $R$ is called \emph{symmetric} if for any  pair $(a,b)$ of
superboxes the following identity holds
\begin{equation}
\label{symm} R(b, a)=R(a, b)^T.
\end{equation}
 Denote by $\Upsilon_i$ the $i$th level of $\Delta$.
Also denote  by $a_i$ and $\sigma_i$  the first and the last
superboxes of the $i$th level $\Upsilon_i$ respectively and by
$r:\Delta\backslash\{\sigma_i\}_{i=1}^d\mapsto \Delta$ the right
shift on the diagram $\Delta$. The last superbox of any level will
be called \emph{special}. For any pair of integers $(i,j)$ such that
$1\leq j<i\leq d$ consider the following tuple of pairs of
superboxes
\begin{equation}
\label{chain}
\begin{split}
~&\bigl (a_j,a_i\bigr),\, \bigl(a_j, r(a_i)\bigr),\, \bigl(r(a_j),
r(a_i)\bigr),\,\bigl(r(a_j), r^2(a_i)\bigr), \ldots,
\bigl(r^{p_i-1}(a_j), r^{p_i-1}(a_i)),
\\
~&\bigl (r^{p_i}(a_j), r^{p_i-1}(a_i)),\ldots, \bigl(r^{p_j-1}(a_j),
r^{p_i-1}(a_i)\bigr).
\end{split}
\end{equation}

Actually the tuple (\ref{chain}) is obtained as follows: the first
pair consists of the last two superboxes of the considered levels,
then until the superbox of the $i$th level will not become special,
each next even pair is obtained from the previous pair of the tuple
by the right shift of the superbox of the $i$th level in the
previous pair and each next odd pair is obtained from the previous
pair of the tuple by the right shift of the superbox of the $j$th
level in the previous pair. When the superbox of the $i$th level
become special, each next pair is obtained from the previous pair of
the tuple by the right shift of the superbox of the $j$th level.

Now we are ready to introduce two crucial notions, which will be
very useful in the formulation of our main Theorem:

\begin{defin}
\label{qnormmap} A symmetric compatible mapping $R:\Delta\times
\Delta\mapsto \rm {Mat}$ is called \underline{quasi-normal}  if the
following two conditions hold: \vskip .1in
\begin{enumerate}
\item
Among all matrices $\mathcal R
(a,b)$, where the superbox $b$ is not higher than the superbox $a$
in the diagram $\Delta$,
the only possible nonzero matrices are the following:
the matrices $\mathcal R
(a,a)$ for
all $a\in \Delta$, the matrices $R
\bigl(a,r(a)\bigr)$, $R
\bigl(r(a),a\bigr)$ for all nonspecial boxes,
and the matrices, corresponding to the pairs, which appear in the
tuples \eqref{chain},
for all $1\leq j<i\leq d$;

\item  The matrix $R\bigl(a,r(a)\bigr)$ is antisymmetric for any nonspecial superbox  $a$.

%
%
%
\end{enumerate}
\end{defin}
\begin{defin}
 \label{normap} A quasi-normal mapping $R:\Delta\times \Delta\mapsto
 \rm {Mat} $ is called \underline{normal} if
it satisfies the
following
condition: for any $1\leq j<i\leq
d$, the matrices, corresponding to the first $(p_j-p_i-1)$ pairs of
the tuple \eqref{chain}, are equal to zero.
\end{defin}

%
%
%

Now let us fix some terminology about the frames in $W$, indexed by
the boxes of the Young diagram $D$.
A frame $\bigl(\{e_\alpha\}_{\alpha\in D},\{f_\alpha\}_{\alpha\in
D}\bigr)$ of $W$ is called \emph{Darboux} or {\it symplectic}, if
for any $\alpha,\beta\in D$ the following relations hold
\begin{equation}
\label{Darboux}
\omega(e_\alpha, e_\beta)=\omega(f_\alpha, f_\beta)=
\omega(f_\alpha, e_\beta)-\delta_{\alpha,\beta}=0,
\end{equation}
where $\delta_{\alpha,\beta}$ is the analogue of the Kronecker index
defined on $D\times D$.
In the sequel it will be  convenient to
divide
a moving frame
$\bigl(\{e_\alpha(t)\}_{\alpha\in D},\{f_\alpha(t)\}_{\alpha\in
D}\bigr)$ of $W$ indexed by the boxes of the Young diagram $D$ into
the tuples of vectors indexed by the supeboxes of the reduction
$\Delta$ of $D$, according to the correspondence between the
superboxes of $\Delta$ and the subcolumns of $D$. More precisely,
given a superbox $a$ in $\Delta$ of size $s$, take all boxes
$\alpha_1, \ldots, \alpha_s$ of the corresponding subcolumn in $D$
in the order from the top to the bottom and denote
\begin{equation*}
E_a(t)=
\bigl(e_{\alpha_1}(t), \ldots, e_{\alpha_s}(t)\bigr),\quad
F_a(t)=
\bigl(f_{\alpha_1}(t),\ldots, f_{\alpha_s}(t)\bigr).
\end{equation*}

In what follows we will suppose that the curve $\Lambda(t)$ is
monotonically  nondecreasing, i.e. the velocity $\dot\Lambda(t)$ is
a nonnegative definite quadratic form for any $t$. The case of
monotonically nonincreasing curve can be treated then by reversing of
time. We restrict ourselves to the monotonic curves just in order to
avoid technicalities both in the formulation and the proof of our
main result (Theorem \ref{maintheor} below). The similar result with
essentially the same proof is valid also for nonmonotonic curves
under additional generic assumptions, which will be introduced in
subsection 3.3 (see Condition (G) there).  In section 4
we point out what changes one has to
make in Theorem \ref{maintheor} in nonmonotonic situation (see
Theorem \ref{maintheor1} below). Note also that Jacobi curves in
sub-Riemannian and, more generally, in sub-Finslerian  geometry are
monotonic, because the corresponding maximized Hamiltonians are
convex on the fibers of $T^*M$ (see the Introduction).
\begin{defin}
\label{normframe} The moving Darboux frame $(\{E_a(t)\}_{a\in
\Delta}, \{F_a(t)\}_{a\in \Delta})$
is called the normal (quasi-normal) moving frame of a monotonically
nondecreasing curve $\Lambda(t)$ with the Young diagram $D$
if
$\Lambda(t)={\rm span}\{E_a(t)\}_{a\in \Delta}$ for any $t$
and there exists an one-parametric family of normal (quasi-normal)
mappings $R_t:\Delta\times\Delta\mapsto \rm {Mat}$ such that the
moving frame $(\{E_a(t)\}_{a\in \Delta}, \{F_a(t)\}_{a\in \Delta})$
satisfies the following structural equation:
\begin{equation}
\label{structeq}
\begin{cases}
 E_a'(t)=E_{l(a)}(t)& \text {if\,\, $a\in
\Delta\backslash\ \mathcal F_1$}\\
E_a'(t)=F_a(t)& \text {if\,\, $a\in \mathcal F_1$}\\
F_a'(t)=\sum\limits_{b\in\Delta}E_b R_t(a,b)-F_{r(a)}& \text {if\,\,
$a\in
\Delta\backslash\ \mathcal S$}\\
F_a'(t)=\sum\limits_{b\in\Delta}E_bR_t(a,b)& \text {if\,\,  $a\in
\mathcal S$}
\end{cases},
\end{equation}
where
$\mathcal F_1$ is the first column of the diagram $\Delta$,
$\mathcal S$ is the set of all its special superboxes, and
$l:\Delta\backslash \mathcal F_1\mapsto\Delta$, $r:\Delta\backslash\
\mathcal S\mapsto \Delta$ are the left and right shifts on the
diagram $\Delta$. The mapping $R_t$, appearing in (\ref{structeq}),
is called the normal (quasi-normal) mapping, associated with the
normal moving frame $(\{E_a(t)\}_{a\in \Delta}, \{F_a(t)\}_{a\in
\Delta})$.
\end{defin}

 With all this terminology we are
ready to formulate our main theorem:

\begin{theor}
\label{maintheor} For any monotonically nondecreasing curve
$\Lambda(t)$ with the Young diagram $D$ in the Lagrange Grassmannian
there exists a normal moving frame
 $(\{E_a(t)\}_{a\in \Delta}, \{F_a(t)\}_{a\in
\Delta})$. A moving frame
$(\{\widetilde E_a(t)\}_{a\in \Delta}, \{\widetilde F_a(t)\}_{a\in
\Delta})$ is a normal moving frame of the curve $\Lambda(\cdot)$ if
and only if for any $1\leq i\leq d$ there exists a constant
orthogonal matrix $U_i$ of size $r_i\times r_i$ such that for all
$t$
\begin{equation}
\label{u} \widetilde E_a(t)= E_a(t)U_i,\quad \widetilde F_a(t)=
F_a(t)U_i, \quad \forall \,a\in \Upsilon_i.
\end{equation}
\end{theor}

Actually, the second statement of this theorem means that if for any
$\bar t$ one collects all possible Darboux frame
$(\{E_a\}_{a\in \Delta}, \{F_a)\}_{a\in \Delta})$ in $W$
such that there exists a normal
moving frame, 
 which coincides with $(\{E_a\}_{a\in \Delta}, \{F_a)\}_{a\in \Delta})$
at $t=\bar t$, then one gets the principle $O(r_1)\times \ldots
\times O(r_d)$ bundle over the curve $\Lambda(t)$ endowed with the
canonical principal connection in the following way: the normal
moving frames are horizontal curves w.r.t. this connection.

\subsection{The canonical splitting and curvature operators}  Before proving Theorem
\ref{maintheor} let us discuss it a little bit.
Take some normal moving frame
$(\{E_a(t)\}_{a\in \Delta}, \{F_a(t)\}_{a\in \Delta})$.
Relations (\ref{u}) imply that for any superbox $a\in\Delta$ of size
$s$ the following $s$-dimensional subspace
\begin{equation}
\label{V} V_a(t)={\rm span}\{E_a(t)\}
\end{equation}
of $\Lambda(t)$ does not depend on the choice of the normal moving
frame. The subspace $V_a$ will be called the \emph
{subspace, associated with the superbox $a$}. 
So, there exists the \emph {canonical splitting of the subspace
$\Lambda(t)$}:
\begin{equation}
\label{cansplit} \Lambda(t)=\bigoplus_{a\in \Delta} V_a(t).
\end{equation}
Moreover, each subspace $V_a(t)$ is endowed with the \emph
{canonical Euclidean structure} such that the tuple of vectors $E_a$
constitute an orthonormal  frame w.r.t. to it. Note that the
canonical splitting is obtained in one of the first steps of the
normalization procedure in the proof of Theorem \ref{maintheor} (see
section 3.4)


Another very important consequence of (\ref{u}) is that the
following subspace
\begin{equation}
\label{trans} \Lambda^{\rm {trans}}(t)=\bigoplus_{a\in\Delta}{\rm
span}\{F_a(t)\}
\end{equation}
does not depend on the choice of the normal moving frame. By
construction,
$W=\Lambda(t)\oplus \Lambda^{\rm {trans}}(t)$
for any $t$. The curve $\Lambda^{\rm {trans}}(t)$ will be called the
\emph{canonical complementary curve of the curve $\Lambda(\cdot)$}.
As we will see in section 5 this notion is crucial for the
construction of the canonical (non-linear) connection for
sub-Riemannian and, more generally, sub-Finsler structures.
\footnote{Note also that this curve is different in general from the
so-called derivative curve $\Lambda^0(\cdot)$, constructed in
\cite{jac1}, which is also intrinsically related to $\Lambda(\cdot)$
such that the  space $\Lambda^0(t)$ is transversal to $\Lambda(t)$
for any $t$. The main disadvantage of the derivative curve
$\Lambda^0(\cdot)$, comparing to the curve $\Lambda^{\rm
{trans}}(\cdot)$, constructed here, is that if one uses it
for the construction of the moving frames intrinsically related to
the curve $\Lambda(\cdot)$ , as was done in \cite{jac1} and
\cite{jac2}, then it is very hard to analyze their  structural
equations and to distinguish a complete system of invariants from it
(in the mentioned papers it was partially done only in the case of
curves of rank 1), while in the present paper we construct the
normal moving frame step by step according to the heuristic rule
that the matrix of its structural equation should be as simple as
possible (should contain as much zeros as possible), which gives the
complete system of invariants automatically.}

Further, we say that a pair $(a,b)$ of superboxes  is
\emph{essential} if $R(a,b)$ is not necessarily zero for a normal
mapping $R:\Delta\times\Delta\mapsto {\rm Mat}$.
Note that this notion depends only on the mutual locations of the
superboxes $a$ and $b$ in the diagram $\Delta$, except the case of
consecutive superboxes $a$ and $b$ in the same level of $\Delta$.
In the last case it depends on the size of the superboxes.  Namely,
the pair $\bigl(a,r(a)\bigr)$ is essential if and only if the size of
$a$ is greater than $1$ (see condition $(1)$ of
Definition \ref{qnormmap}).

Assume that  $R_t :\Delta\times\Delta\mapsto {\rm Mat}$ and
$\widetilde R_t :\Delta\times\Delta\mapsto {\rm Mat}$ are the normal
mappings, associated with normal moving frames $(\{E_a(t)\}_{a\in
\Delta}, \{F_a(t)\}_{a\in \Delta})$ and $(\{\widetilde
E_a(t)\}_{a\in \Delta}, \{\widetilde F_a(t)\}_{a\in \Delta})$,
which are related by (\ref{u}).
Then from (\ref{structeq}) and (\ref{u}) it follows immediately
that
\begin{equation}
\label{Rtrans} \widetilde R_t(a,b)=U_j^{-1} R_t(a,b)U_i, \quad a\in
\Upsilon_i, b\in \Upsilon_j.
\end{equation}
 The last relation means actually that for any essential
 pair $(a,b)$ of superboxes the
linear mapping ${\mathfrak R}_t (a,b): V_a\mapsto V_b$, having the
matrix $R_t(a,b)$ w.r.t. the bases $E_a$ and $E_b$ of $V_a$ and
$V_b$ respectively, does not depend on the choice of a normal moving
frame.\footnote {Here we restrict ourselves to essential pairs,
because for nonessential pairs such linear mappings are zeros
automatically.} The linear mapping ${\mathfrak R}_t (a,b)$ will be
called the \emph {$(a,b)$-curvature mapping of the curve
$\Lambda(\cdot)$}.

The only nontrivial blocks in the matrix of the structural equations
for the normal moving frames correspond to $(a,b)$-curvature
mappings. Hence the tuple of all $(a,b)$-curvature mappings
constitute a kind of complete system of symplectic invariants of the
curve. For precise formulation of this statement it is convenient to
use the notion of quivers and their representations (\cite{quiv}).
Recall that a quiver is an oriented graph, where loops and multiple
arrows between two vertices are allowed. A representation of a
quiver assigns a vector space  $X_\alpha$  to each vertex $\alpha$
of the quiver and a linear mapping  from $X_\alpha$ to $X_\beta$ to
each arrow of the quiver, connecting a vertex $\alpha$ with  a
vertex $\beta$.

Take the quiver $\mathfrak Q_D$ such that its vertices are levels of
the diagram $\Delta$ and the set of arrows from the level
$\Upsilon_i$ to the level $\Upsilon_j$ is parameterized by essential
pairs $(a,b)\in \Upsilon_i\times \Upsilon_j$. A representation of
the quiver $\mathfrak Q_D$ will be called \emph{compatible  with the
Young diagram $D$} if for any $1\leq i\leq d$ the space of the
representation corresponding to the vertex $\Upsilon_i$ is a
$r_i$-dimensional Euclidean space and the linear mappings $\mathcal
R(a,b)$ of the representation corresponding to the arrows $(a,b)$
satisfy the following relations: $\mathcal R(a,b)^*=\mathcal R(b,a)$
and $\mathcal R\bigl(a,r(a)\bigr)$ are antisymmetric w.r.t. the
corresponding Euclidean structure.

The subspaces $V_a(t)$ for any $t$ and any $a\in \Upsilon_i$ are
naturally identified together with the canonical Euclidean structure
on them ($V_{a_1}(t_1)\sim V_{a_2}(t_2)$ by sending $E_{a_1}(t_1)$
to $E_{a_2}(t_2)$). Therefore, we can identify all these spaces with one
Euclidean space, which will be denoted by $\mathcal X_i$.
The tuple of spaces  $\mathcal X_i$ and the $(a,b)$-curvatures
mappings of the curve $\Lambda(t)$, considered as elements of ${\rm
Hom }(\mathcal X_i, \mathcal X_j)$ for $(a,b)\in\Upsilon_ i\times
\Upsilon_j$,  define the one-parametric family $\mathfrak R_t$ of
compatible representations of the quiver $\mathfrak Q_D$. This
family will be called the \emph {quiver of curvatures of  the curve
$\Lambda(t)$}. Here the linear mappings corresponding to the arrows
of the quiver depend on $t$, while the linear spaces, corresponding
to its vertices, are independent of $t$. In the sequel we will
consider only this type of one-parametric families of
representations of quivers. Two families $\Xi_1(t)$ and $\Xi_2(t)$
of compatible representations of the quiver $\mathfrak Q_D$ are
called isomorphic, if there exists a tuple of isometries
(independent of $t$) between the corresponding spaces of the
representations, conjugating all corresponding linear mappings.
If the sizes of all superboxes in $\Delta$ are equal to 1, then the
normal moving frames of the curve are defined up to the discrete
group ($U_i$ in \eqref{u} are scalars, which are equal to $1$ or
$-1$) and all $(a,b)$-curvature mappings are determined by scalar
functions of $t$, which are symplectic invariants of the curve.
These scalar functions will be called, for short,
\emph{$(a,b)$-curvatures}. Besides, the compatible representations of
the quiver $\mathfrak Q_D$ is in one-to-one correspondence with
tuples of numbers parameterized by the essential pairs of $\Delta$
(which is equal to $D$ in the considered case). The following
theorem is the direct consequence of the structural equations for
normal moving frames and Theorem \ref{maintheor}:

\begin{theor}
\label{size1} For the given one-parametric family $\Xi(t)$ of
representations of the quiver $\mathfrak Q_D$ compatible with the
Young diagram $D$ with $|D|$ boxes there exists the unique, up to a
symplectic transformation, monotonically nondecreasing curve
$\Lambda(t)$ in the Lagrange Grassmannian of $2|D|$-dimensional
symplectic space with the Young diagram $D$ such that the quiver of
curvatures of
$\Lambda(t)$ is isomorphic to $\Xi(t)$.
If, in addition,  all rows of
$D$ have different length,
then given a tuple of smooth functions
$\{\rho_{a,b}(t):(a,b)\in\Delta\times \Delta, (a,b) \,\,\text{is an
essential pair} \}$ there exists the unique, up to a symplectic
transformation, monotonically nondecreasing curve $\Lambda(t)$ in
the Lagrange Grassmannian of $2|D|$-dimensional symplectic space
with the Young diagram $D$ such that for any essential pair
$(a,b)\in\Delta\times \Delta$ and any $t$ its $(a,b)$-curvature at
$t$ coincides with $\rho_{a,b}(t)$.
\end{theor}

Finally note  that rank 1 curves in Lagrange Grassmannians,
considered in \cite{zel}, have the Young diagrams, consisting of
just one row, and the main results of the mentioned paper (Theorems
2 and 3 there) are very particular cases of Theorems \ref{maintheor}
and \ref{size1} here. In this case the pair $(a,b)$ of superboxes is
essential if and only if $a=b$.
%
%

\section{Proof of Theorem 1}
\setcounter{equation}{0}

The proof consists of several steps.
\subsection{Contractions of the curve $\Lambda(\cdot)$}
We start with some general constructions for curves in
Grassmannians. Given a curve $\Lambda(\cdot) $ in the Grassmannian
$G_k(W)$, for any $\tau$ we will construct a monotonic sequence of
subspaces of $\Lambda(\tau)$ in addition to the extensions
$\Lambda^{(i)}$. For this let $\Lambda_{(0)}(t)=\Lambda(t)$ and
recursively
\begin{equation}
\label{primeicontr} \Lambda_{(i)}(\tau)=
\left\{v\in  \Lambda_{(i-1)}(\tau):
\begin{array}{l}
\exists\,\ell\in\mathfrak S(\Lambda_{(i-1)})\,\,
\text {with}\,\, \ell(\lambda)=v\\
\text {such that}\,\,
\ell'\bigl(\tau)\in \Lambda_{(i-1)}(\tau)
\end{array}\right\}
\end{equation}
where, by analogy with above, $\mathfrak S(\Lambda_{(i)})$, $i\geq
0$, is the set of all smooth curves $\ell(t)$ in $W$ such that
$\ell(t)\in \Lambda_{(i-1)}(t)$ for any  $t$.
The subspaces $\Lambda_{(i)}(\tau)$ are called {\it the $i$th
contraction} of the curve $\Lambda(\cdot)$ at the point $\tau$.
Under the identification $T_\Lambda G_k(W)\sim {\rm
Hom}\, (\Lambda,W/\Lambda)$ the first contraction
$\Lambda_{(1)}(\tau)$ is exactly the kernel of the velocity
$\dot\Lambda(\tau)$, $
\Lambda_{(1)}(\tau) ={\rm Ker}\,\dot\Lambda(\tau)$.
In particular, it implies that
\begin{equation}
\label{kerim} \dim \Lambda^{(1)}(\tau)- \dim \Lambda(\tau)= \dim
\Lambda(\tau)-\dim \Lambda_{(1)}(\tau).
\end{equation}
Indeed, the righthand side of (\ref{kerim}) is equal to $\dim
\bigl({\rm Im} \,\dot\Lambda(\tau)\bigr)$, while the lefthand side
is equal to $\dim \Lambda(\tau)-\dim \bigl({\rm
Ker}\,\dot\Lambda(\tau)\bigr)$.

Note also that in \eqref{primeicontr} one can replace the quantor
$\exists$ by $\forall$, because the existence of a curve $\ell\in
\mathfrak S(\Lambda_{(i-1)}) $  with $\ell(\tau)=v$ and $
\ell^\prime(\tau)\in \Lambda _{(i-1)}(\tau)$ implies that any smooth
curve $\tilde\ell\in\mathfrak S(\Lambda_{(i-1)})$ with
$\tilde\ell(\tau)=v$ satisfies $\tilde \ell^\prime(\tau)\in \Lambda
_{(i-1)}(\tau)$.
Note that the following relations follow directly from the
definitions
\begin{equation}
\label{propbas}
\bigl(\Lambda_{(i)}(\tau)\bigr)_{(1)}=\Lambda_{(i+1)}(\tau),\quad
\bigl(\Lambda_{(i)}(\tau)\bigr)^{(1)}\subseteq \Lambda_{(i-1)}(\tau)
\end{equation}

If we suppose that $\Lambda(\cdot)$ is a curve in Lagrange
Grassmannian of the symplectic space  $W$, then
the symplectic structure gives an additional relation between the
$i$th extension and the $i$th contraction. Namely, given a subspace
$L\subset W$ denote by $L^\angle$ its skew-symmetric complement,
i.e. $L^\angle=\{v\in W: \omega (v,l)=0\,\, \forall l\in L\}$.
\begin{lemma}
\label{subsuplem} The subspaces $\Lambda_{(i)}(\tau)$ is a
skew-symmetric complement of the subspace $\Lambda^{(i)}(\tau)$ for
any $\tau$, namely
\begin{equation}
\label{subsup}
\Lambda_{(i)}(\tau)=\Bigl(\Lambda^{(i)}(\tau)\Bigr)^\angle, \quad
\forall \tau.
\end{equation}
\end{lemma}
{\bf Proof.}
 We proceed the proof by induction on $i$. For $i=0$
there is nothing to prove, because $\Lambda(\tau)$
($=\Lambda^{0}(\tau)=\Lambda_{0}(\tau)$ by definition) is a
Lagrangian subspace. Assume that (\ref{subsup}) is valid for $i=\bar
i-1$ and prove it for $i=\bar i$, $\bar i\geq 1$. Indeed, if $v\in
\Lambda_{(\bar i)}(\tau)$, then by definition there exists a regular
curve of vectors $v(t)$ such that $v(t) \in \Lambda_{(\bar i-1)}(t)$
for any $t$  close to $\tau$, $v(\tau)=v$ and $v'(\tau)\in \Lambda_{
(\bar i-1)}(\tau)$. Let us prove that $v\in \bigl(\Lambda_{(\bar
i)}(\tau)\bigr)^\angle$. For this take $v_1\in \Lambda^{(\bar
i)}(\tau)$. Then by definition there exist a curve of vectors $w(t)$
in $W$ such that $w(t)\in \Lambda^{(\bar i-1)}(t)$ for any $t$ close
to $\tau$ and  $w'(\tau)=v_1$. By induction hypothesis
$\omega\bigl(v(t), w(t)\bigr)=0$. Differentiating the last identity
at $t=\tau$ we get
\begin{equation}
\label{prom1} \omega (v, v_1)=-\omega\bigl(v'(\tau),w(\tau)\bigr)=0.
\end{equation}
(the last equality holds because of the relations $v'(\tau)\in
\Lambda_{ (\bar i-1)}(\tau)$, $w(\tau)\in \Lambda^{(\bar
i-1)}(\tau)$ and the induction hypothesis). Since \eqref{prom1}
holds for any $v_1\in\Lambda^{(\bar i)}(\tau)$, we get that $v\in
\bigl(\Lambda_{(\bar i)}(\tau)\bigr)^\angle$. So, we have proved
that
$\Lambda_{(i)}(\tau)\subset\Bigl(\Lambda^{(i)}(\tau)\Bigr)^\angle$.

Now let us prove the inclusion in the opposite direction. Suppose
that $v\in\Bigl(\Lambda^{(\bar i)}(\tau)\Bigr)^\angle$.
Take any $w\in \Lambda^{(\bar i-1)}(\tau)$ and a curve of vectors
$w(t)$ in $W$ such that $w(t)\in \Lambda^{(\bar i-1)}(t)$ for any
$t$ close to $\tau$ and $w(\tau)=w$. Then by definition $w'(\tau)\in
\Lambda^{(\bar i)}(\tau)$ and by our assumptions
\begin{equation}
\label{prom2} \omega \bigl(v,w'(\tau)\bigr)=0.
\end{equation}
On the other hand, since $\Lambda^{(\bar i-1)}(\tau)\subset
\Lambda^{(\bar i)}(\tau)$, then $\Bigl(\Lambda^{(\bar
i)}(\tau)\Bigr)^\angle\subset \Bigl(\Lambda^{(\bar
i-1)}(\tau)\Bigr)^\angle=\Lambda_{(\bar i-1)}(\tau)$ (the last
equality is our induction hypothesis). So, $v\in\Lambda_{(\bar
i-1)}(\tau)$. Take a curve of vectors $v(t)$ in $W$  such that
$v(t)\in \Lambda_{(\bar i-1)}(t)$ for any $t$ close to $\tau$ and
$v(\tau)=v$. Then by induction hypothesis $\omega(v(t),w(t))=0$ for
any $t$ close to $\tau$. Differentiating the last identity at
$t=\tau$ and using (\ref{prom2}) we get that $\omega\bigl(v'(\tau),
w\bigr)=0$. Since the last identity holds for any $w\in
\Lambda^{(\bar i-1)}(\tau)$, then $v'(\tau)\in \Bigl(\Lambda^{(\bar
i-1)}(\tau)\Bigr)^\angle =\Lambda_{(\bar i-1)}(\tau)$ (the last
equality is our induction hypothesis). So, $v\in \Lambda_{(\bar
i)}(\tau)$, which implies the inclusion $\Bigl(\Lambda^{(\bar
i)}(\tau)\Bigr)^\angle\subset\Lambda_{(\bar i)}(\tau)$. The proof of
the lemma is completed.
$\Box$

\subsection{
Filling the Young diagram $D$ by bases of $\Lambda(t)$}
 As before, assume that the reduced diagram
$\Delta$ of the curve consists of $d$ level,  the number of
superboxes in the $i$th level of the diagram $\Delta$ is equal to
$p_i$,  and their sizes are equal to $r_i$.
By our assumptions $\Lambda^{(p_1)}(t)=W$, which together with
(\ref{subsup}) implies that
\begin{equation}
\label{esig00} \Lambda_{(p_1)}(t)=0,\quad \dim
\Lambda_{(p_1-1)}(t)=r_1.
\end{equation}
Denote also by $\sigma_i$ the special (i.e. the last) superbox of
the $i$th level of $\Delta$.
From the second relation of
(\ref{propbas}) it follows that
\begin{equation}
\label{esig01}
\Bigl(\Lambda_{(p_{i})}\Bigr)^{(1)}(t)
\subseteq\Lambda_{(p_{{i}}-1)}(t),\quad \forall \,1\leq i \leq q
\end{equation}


For any $1\leq i\leq d$ choose  a complement $\widetilde
V_{\sigma_i}(t)$ of the subspace
$\Bigl(\Lambda_{(p_{i})}\Bigr)^{(1)}(t)$ in the space
$\Lambda_{(p_i-1)}(t)$ (smoothly w.r.t. $t$):
\begin{equation}
\label{esig02}
\Lambda_{(p_{i}-1)}=
\Bigl(\Lambda_{(p_{i})}\Bigr)^{(1)}(t) \oplus \widetilde
V_{\sigma_i}(t).
\end{equation}
Note that from (\ref{esig00})
it follows
that $\widetilde V_{\sigma_1}(t)=\Lambda_{(p_1-1)}(t)$. Let
$\widetilde\Delta$ be the diagram, obtained from $\Delta$ by joining
to $\Delta$ one more column from the left, having the same length as
the first column of $\Delta$. The boxes of $\widetilde \Delta$ will
be called superboxes as well.  For any $1\leq i\leq d$ take a tuple
of vectors $E_{\sigma_i}(t)$, constituting a basis of $\widetilde
V_{\sigma_i}(t)$ (smoothly in $t$). Then to any superbox of $\tilde
\Delta$ we will assign a tuple of vectors in the following way
\begin{equation}
\label{esig6}
E_{l^{^{j}}(\sigma_i)}(t)\stackrel{def}{=}E_{\sigma_i}^{(j)}(t),\quad
\forall\, 0\leq j\leq p_i,
\end{equation}
where $l$ is the left shift on the diagram $\widetilde\Delta$.
\begin{lemma}
\label{filling} 
Assume that a superbox $a\in\widetilde \Delta$ lies in the $j(a)$th
column and $i(a)$th level of the diagram $\widetilde\Delta$ and let
${\rm Ov}_a$ be the set of all superboxes, lying over $a$ in the
column of $a$. Then the following relations hold
\begin{equation}
\label{esigm}
\begin{split}
~&
\bigl\{E_a(t)\bigr\}\bigcap\Bigl(\Bigl(\bigoplus_{b\in {\rm Ov}_a}
{\rm
span}\bigl\{E_b(t)\bigr\}\Bigr)\oplus\Lambda_{(j(a)-1)}(t)\Bigr)=0
,\\
~&\dim\, {\rm span} \bigl\{E_a(t)\bigr\}=\dim\, {\rm span}
\bigl\{E_{\sigma_{i(a)}}(t)\bigr\}=r_{i(a)}.
\end{split}
\end{equation}
\end{lemma}

{\bf Proof}.
 Let $\prec$ be the order on the set of superboxes of
the diagram $\widetilde \Delta$, defined as follows: $b_1\prec b_2$
if either $b_1$ is higher than $b_2$ in $\widetilde \Delta$ or they
are on the same level, but $b_1$ is located from the right to $b_2$
(or, equivalently, either $i(b_1)<i(b_2)$ or $i(b_1)=i(b_2)$, but
$j(b_1)>j(b_2)$). Let us prove (\ref{esigm}) by induction on the set
of superboxes of  the diagram $\widetilde \Delta$ with the
introduced order $\prec$. For $a=\sigma_1$ relations (\ref{esigm})
follow immediately from (\ref{esig00}). Now assume that
(\ref{esigm}) is true for any superbox $a\in\widetilde D$ such that
$a\prec \mathfrak \sigma$ and prove it for $a=\sigma$.
We have the following two cases:

{\bf 1. The superbox $\sigma$ is special.} In this case by induction
hypothesis it is easy to show that
\begin{equation}
\label{esig03} \Bigl(\bigoplus_{b\in {\rm Ov}_{\sigma}} {\rm
span}\bigl\{E_{b}(t)\bigr\}\Bigr)\oplus\Lambda_{(p_{i(\sigma)})}(t)=
\Bigl(\Lambda_{(p_{i(\sigma)})}\Bigr)^{(1)}(t)
\end{equation}
This together with (\ref{esig02}) and the definitions of the numbers
$r_i$ implies (\ref{esigm}) for $a=\sigma$.

{\bf 2. The superbox $\sigma$ is not special.} Using our induction
assumptions we can choose a subspace $C(t)$ of
$\Lambda_{(j(\sigma)-1)}(t)$ smoothly w.r.t. $t$ such that

\begin{equation}
\label{esig3} \Lambda_{(j(\sigma)-1)}(t)=\Bigl(\bigoplus_{b\in {\rm
Ov}_{r(\sigma)}}{\rm span}\bigl\{E_b(t)\bigr\}\Bigr)\oplus{\rm
span}\bigl\{E_{r(\sigma)}(t)\bigr\}
\oplus\Lambda_{(j(\sigma))}(t)\oplus C(t),
\end{equation}
where as before $r(\sigma)$ is the superbox, located from the right
to $\sigma$ in $\widetilde\Delta$.

From (\ref{kerim}), the first relation of (\ref{propbas}), and
(\ref{esig3}) it follows that
\begin{equation}
\label{esig4}
\begin{split}
~&\dim (\Lambda_{(j(\sigma)-1)})^{(1)}(t)-\dim
\Lambda_{(j(\sigma)-1)}(t)= \dim\Lambda_{(j(\sigma)-1)}(t)-\dim
\Lambda_{(j(\sigma))}(t)=\\~&\sum_{b\in {\rm Ov}_{r(\sigma)}\cup
r(\sigma)}\dim {\rm span}\bigl\{E_b(t)\bigr\}+ \dim
C(t)=\sum_{k=1}^{i(\sigma)}r_k+\dim C(t).
\end{split}
\end{equation}
On the other hand, using definitions (\ref{primeiext}),
(\ref{primeicontr}), (\ref{esig6}), \eqref{esig3}, the induction
hypothesis,
and relation (\ref{kerim}) one gets easily that
\begin{equation}
\label{esig5}
\begin{split}
~&\dim (\Lambda_{(j(\sigma)-1)})^{(1)}(t)-\dim
\Lambda_{(j(\sigma)-1)}(t)\leq \sum_{k=1}^{i(\sigma)-1}r_k+
\left(\dim {\rm span}\bigl\{E_{r(\sigma)}(t),
E_\sigma(t)\bigr\}-\right.\\
~&\left. \dim {\rm span}\bigl\{E_{r(\sigma)}(t)\bigr\}\right)+
\bigl(\dim C^{(1)}(t)-\dim C(t)\bigr)\leq
\sum_{k=1}^{i(\sigma)}r_k + \dim C(t).
\end{split}
\end{equation}
If for $a=\sigma$ one of the identities in (\ref{esigm}) does not
hold, then in the chain of the inequalities (\ref{esig5}) there is
at least one strong inequality, which is in the contradiction with
(\ref{esig4}). So, the identities (\ref{esigm}) are valid for
$a=\sigma$, which completes the proof of (\ref{esigm}) by induction.
$\Box$
\medskip

 Let $\mathcal F_k$ be the $k$th column of the diagram
$\Delta$. From Lemma \ref{filling} it follows easily the following
 \begin{cor}
 The following splittings hold for any $0\leq j\leq p_1$
\begin{equation}
 \label{splitfl}
\Lambda_{(j)}(t)=\bigoplus_{a\in \bigcup_{s=j+1}^{p_1} \mathcal
F_{s}} {\rm span} \{E_a(t)\}, \quad
\bigl(\Lambda_{(j)}\bigr)^{(1)}(t)=\bigoplus_{a\in
\bigcup_{s=j+1}^{p_1} \mathcal F_s\cup l(F_{j+1})} {\rm span}
\{E_a(t)\}
\end{equation}
In particular, $\Lambda(t)=\bigoplus_{a\in  \Delta} {\rm span}
\{E_a(t)\}$.
\end{cor}
One can imagine that we fill the diagram $\Delta$ (or the original
diagram $D$) by columns $E_a(t)^T$ by choosing bases of the
subspaces $\widetilde V_{\sigma_i}$, satisfying \eqref {esig02}, and
by differentiating these bases as in \eqref{esig6}. Tuples
$\{E_a(t)\}_{a\in\Delta}$, obtained in this way, will be called
\emph{fillings of the Young diagram D, associated  with the curve
$\Lambda(\cdot)$}. The flag $0={\Lambda_{(p_1}}(t)\subset
{\Lambda_{(p_1-1}}(t) \ldots\subset \Lambda_{(0)}(t)=\Lambda(t)$ can
be recovered from this filling by the first relation of
(\ref{splitfl}). In particular, this flag (and therefore the curve
$\Lambda(\cdot)$ itself) can be recovered
from the curves $t\mapsto V_{\sigma_i}(t)$, $1\leq i\leq d$ by
taking the corresponding extensions of them.

\subsection{The canonical complement of
$\Bigl(\Lambda_{(p_{i})}\Bigr)^{(1)}(t)$ in $\Lambda_{(p_{i}-1)}(t)$
and the canonical Euclidean structure on it} In the previous
subsection we took some complements $\widetilde V_{\sigma_i}$ of the
subspaces $\Bigl(\Lambda_{(p_{i})}\Bigr)^{(1)}(t)$
in the spaces $\Lambda_{(p_{i}-1)}(t)$. In the present section we
will show that such complements can be chosen canonically for a
curve $\Lambda(t)$ with the Young diagram $D$, satisfying the
following additional assumption:
\medskip

{\bf Condition (G)} \emph {For any $1\leq i\leq d-1$ and any $t$ the
rank of the restriction of the quadratic form $\dot \Lambda(t)$ to
the subspace $\bigl(\Lambda_{(p_i-1)}\bigr)^{(p_i-1)}(t)$ is equal
to $\displaystyle{\sum_{k=1}^i r_k}$,}
\begin{equation}
\label{condG} \forall 1\leq i\leq d-1 \text { and } \forall\,t:\quad
{\rm
rank}\,\left(\dot\Lambda(t)|_{(\Lambda_{(p_i-1)})^{(p_i-1)}(t)}\right)=
\sum_{k=1}^i r_k.
\end{equation}
\medskip
Since ${\rm Ker}\,\dot\Lambda(t)=\Lambda_{(1)}(t)$ and
$\Lambda_{(p_i-1)})^{(p_i-2)}(t)\subset \Lambda_{(1)}(t)$ (as a
consequence of \eqref{propbas}), any curve $\Lambda(t)$ with the
Young diagram $D$ satisfies: ${\rm
rank}\,\left(\dot\Lambda(t)|_{(\Lambda_{(p_i-1)})^{(p_i-1)}(t)}\right)\leq
\sum_{k=1}^i r_k$ for any $1\leq i\leq d$. It implies easily that
germs of curves, satisfying condition (G), are generic among all
germs of curves with given Young diagram $D$. Besides, it is clear
that curves with rectangular Young diagram satisfy condition (G)
automatically (condition (G) is empty in this case).

\begin{lemma}
\label{condGlem} Any monotonic curve $\Lambda(t)$ with the Young
diagram $D$ satisfies condition (G).
\end{lemma}

{\bf  Proof}. For definiteness, assume that the curve  $\Lambda(t)$
is monotonically nondecreasing. Take a filling
$\{E_a(t)\}_{a\in\Delta}$ of the Young diagram $D$, associated with
the curve $\Lambda(\cdot)$. Let
\begin{equation}
\label{Z} Z_i(t)={\rm span}\{ E_{\sigma_k}^{(p_k-1)}(t) \} _{k=1}^i,
\quad
1\leq i\leq q.
\end{equation}
 It is clear that $\{Z_i(t)\}_{i=1}^d$ is a
monotonically increasing (by inclusion) sequence of subspaces for any $t$. As a
consequence of Lemma \ref{filling}, we have
\begin{eqnarray}
&~& \dim\, Z_i(t)=\sum_{k=1}^i r_k,\label{mon0}\\
&~& \label{mon1} (\Lambda_{(p_i-1)})^{(p_i-1)}(t)=
\Bigl((\Lambda_{(p_i-1)})^{(p_i-1)}(t)\cap
\Lambda_{(1)}(t)\Bigr)\oplus Z_i(t)
\end{eqnarray}
Since ${\rm Ker}\,\dot\Lambda(t)=\Lambda_{(1)}(t)$, we get from
\eqref{mon1} that
\begin{equation}
\label{mon2} {\rm
rank}\,\left(\dot\Lambda(t)|_{(\Lambda_{(p_i-1)})^{(p_i-1)}(t)}\right)={\rm
rank}\,\left(\dot\Lambda(t)|_{Z_i(t)}\right).
\end{equation}
Besides, from monotonicity the quadratic form
$\dot\Lambda(t)|_{Z_d(t)}$ is positive definite. Hence, the
quadratic forms $\dot\Lambda(t)|_{Z_i(t)}$ are positive definite as
well. Then the lemma follows form  relations \eqref{mon0} and
\eqref{mon2}.
$\Box$
\medskip

Now define the following subspaces of the ambient symplectic space
$W$:
\begin{equation}
\label{sympi}
W_i(t)=(\Lambda_{(p_1-1)}(t))^{(2p_1-1)}+(\Lambda_{(p_2-1)}(t))^{(2p_2-1)}+\ldots+(\Lambda_{(p_i-1)}(t))^{(2p_i-1)}.
\end{equation}

\begin{lemma}
\label{simpilem} If a curve $\Lambda(t)$ with the Young diagram $D$
satisfies condition (G), then for any $1\leq i\leq d$ the
restriction of the symplectic form $\omega$ to the subspace $W_i(t)$
is nondegenerated and $dim W_i=2\displaystyle{\sum _{k=1}^i p_k
r_k}$.
\end{lemma}
{\bf Proof}. The proof of the lemma is by induction w.r.t. $i$.
First let us introduce some notations. Let $\overline\Delta$ be the
diagram obtained from $\Delta$ by the reflection w.r.t. its left
edge. We will work with the diagram $\Delta\cup\overline\Delta$,
which is symmetric w.r.t. the left edge of the diagram $\Delta$.
Similar to
above, we will denote by $l$
the left 
shift on the diagram $\Delta\cup\overline\Delta$. If $S$ is a subset
of the diagram $\Delta$,   we will denote by $\bar S$ the subset of
$\bar\Delta$, obtained by the reflection of $S$ w.r.t. the left edge
of $\Delta$. Also in the sequel, given two tuples of vectors
$V_1=(v_{11}, \ldots, v_{1n_1})$ and $V_2=(v_{21}, \ldots,
v_{2n_2})$ by $\omega(V_1, V_2)$ we will mean the $n_1\times
n_2$-matrix with the $(i,j)$-entry equal to $\omega(v_{1i},
v_{2j})$. Take a filling $\{E_a(t)\}_{a\in\Delta}$ of the Young
diagram $D$, associated with the curve $\Lambda(\cdot)$. Define
tuples $E_a$ also for $a\in\bar \Delta$ in the following way:
$E_{l^j(a_i)}= E_{a_i}^{(j)}(t)$, $1\leq j\leq p_i$,  where, as
before,  $a_i$ is the first superbox in the $i$th level $\Upsilon_i$
of $\Delta$. By definition $W_i(t)={\rm span}\{
E_a(t)\}_{a\in\cup_{k=1}^i\Upsilon_k\cup\bar\Upsilon_k}$.

{\bf 1.} Let us prove the lemma for $i=1$. By condition (G) the
matrix $\omega\bigl(E_{\sigma_1}^{(p_1)}(t),
E_{\sigma_1}^{(p_1-1)}(t)\bigr)$ is nonsingular. On the other hand,
since $\Lambda_{(1)}(t)=\Bigl(\Lambda^{(1)}(t)\Bigr)^\angle$, one
has $\omega\bigl(E_{\sigma_1}^{(p_1)}(t),
E_{\sigma_1}^{(p_1-2)}(t)\bigr)\equiv 0$. Differentiating the last
identity, we get $\omega\bigl(E_{\sigma_1}^{(p_1+1)}(t),
E_{\sigma_1}^{(p_1-2)}(t)\bigr)=-\omega\bigl(E_{\sigma_1}^{(p_1)}(t),
E_{\sigma_1}^{(p_1-1)}(t)\bigr)$. In the same way, using
\eqref{subsup}, it is easy to obtain that
\begin{equation*}
\omega\bigl(E_{\sigma_1}^{(p_1+i)}(t),
E_{\sigma_1}^{(p_1-i-1)}(t)\bigr)=(-1)^i\omega\bigl(E_{\sigma_1}^{(p_1)}(t),
E_{\sigma_1}^{(p_1-1)}(t)\bigr).
\end{equation*}
In particular, all matrices $\omega\bigl(E_{\sigma_1}^{(p_1+i)}(t),
E_{\sigma_1}^{(p_1-i-1)}(t)\bigr)$ are nonsingular. Therefore the
matrix with the entries, which are equal to the value of the form
$\omega$ on all pairs of vectors from the tuple $\{
E_a(t)\}_{a\in\Upsilon_1\cup\bar\Upsilon_1}$, is block-triangular
w.r.t. the nonprincipal diagonal with nonsingular blocks on the
nonprincipal diagonal. This implies that the tuple $\{
E_a(t)\}_{a\in\Upsilon_1\cup\bar\Upsilon_1}$ constitutes the basis
of $W_1$ and the form $\omega|_{W_1}$ is nondegenerated, which
completes the proof of the statement of the lemma in the case $i=1$.

{\bf 2.} Now assume that the statement of the lemma holds for
$i=i_0-1$ and prove it for $i=i_0$. Let $\Delta_i$ be the subdiagram
of $\Delta$, consisting of the first $i$ rows of $\Delta$,
$\Delta_i=\displaystyle{\bigcup_{k=1}^i \Upsilon_k}$. Divide the
diagram $\Delta_{i_0}\cup\bar\Delta_{i_0}$ on four parts
$\{A_k\}_{k=1}^4$: $A_1$ is a union of the last $p_1-p_{i_0}$
columns of the diagram $\Delta_{i_0}$, $A_2$ is obtained by the
reflection of $A_1$ w.r.t. the left edge of $\Delta_{i_0}$, i.e.
$A_2=\bar A_1$, $A_3=\Delta_{i_0-1}\backslash (A_1\cup A_2)$, and
$A_4=\Upsilon_{i_0}$.

Set $C_k(t)={\rm span}\{ E_a(t)\}_{a\in A_k}$, $k=1,\ldots, 4$. Note
that from \eqref{splitfl} it follows that
$C_1(t)=\Lambda_{(p_{i_0})}(t)$. By constructions
$W_{i_0}(t)=C_1(t)+C_2(t)+C_3(t)+C_4(t)$ and
$W_{i_0-1}=C_1(t)+C_2(t)+C_3(t)$. Moreover, by induction hypothesis
\begin{eqnarray}
&~& W_{i_0-1}(t)=C_1(t)\oplus C_2(t)\oplus C_3(t)
\label{omeg2}\\&~&\label{omeg3} C_1(t)^\angle\cap
W_{i_0-1}(t)=C_1(t)\oplus C_3(t).
\end{eqnarray}
The last two identities follow just from comparison of dimensions.
Besides, using \eqref{subsup}, one has also that
\begin{equation}
\label{omeg4} C_1+C_3+C_4\subset C_1(t)^\angle.
\end{equation}
%
Assume that $x\in {\rm Ker} \,\omega|_{W_{i_0}(t)}$,
$\displaystyle{x=\sum_{k=1}^4 x_k}$, where  $x_k\in C_k(t)$. Then
\eqref{omeg4} implies that $\omega(v,x)=\omega (v,x_2)=0$ for any
$v\in C_1(t)$. This together with  \eqref{omeg2} and \eqref{omeg3}
yields that $x_2=0$.

Further, by the same arguments as in the proof of the case $i=1$,
applied for the tuple $\{E_a\}_{a\in \mathcal F_{p_{i_0}}\cap
\Delta_{i_0}}$ instead of the tuple $E_{\sigma_1}$, one obtains from
\eqref{condG} for $i=i_0$ that $\omega|_{C_3(t)+C_4(t)}$ is
nondegenerated and $\dim
(C_3(t)+C_4(t))=2p_{i_0}\sum_{k=1}^{i_0}r_k$. The latter implies
that $C_1(t)\cap C_3(t)=0$. Besides,  from \eqref{omeg4} it follows
that $\omega(v,x)=\omega (v,x_3+x_4)=0$ for any $v\in
C_3(t)+C_4(t)$, which together with two previous sentences implies
that $x_3=x_4=0$. Therefore  $x\in C_1(t)\subset W_{i_0-1}(t)$,
which implies that $x=x_1=0$ by induction hypothesis. This yields
that the form $\omega|_{W_{i_0}(t)}$ is nondegenerated. Moreover,
from the same arguments it follows that the condition
$\displaystyle{\sum_{k=1}^4 x_k=0}$ implies that $x_k=0$ for any
$1\leq k\leq 4$. Hence $W_{i_0}(t)=C_1(t)\oplus C_2(t)\oplus
C_3(t)\oplus C_4(t)$ and the statement of the lemma about the
dimension of $W_{i_0}(t)$ holds. The proof of the lemma is
completed.
$\Box$
\medskip

Finally, let
\begin{equation}
\label{cancomp} V_i(t)=\Lambda_{(p_i-1)}(t)\cap W_{i-1}(t)^\angle
\end{equation}
As a direct consequence of Lemma \ref{simpilem}, we get that the
subspace $V_i(t)$ is complementary to
$\Bigl(\Lambda_{(p_{i})}\Bigr)^{(1)}(t)$ in
$\Lambda_{(p_{i}-1)}(t)$,
\begin{equation}
\label{cancomp1}
\Lambda_{(p_{i}-1)}(t)=\Bigl(\Lambda_{(p_{i})}\Bigr)^{(1)}(t)\oplus
V_i(t). \end{equation}
The subspaces $V_i(t)$, defined by \eqref{cancomp}
will be called
the \emph {canonical complement of
$\Bigl(\Lambda_{(p_{i})}\Bigr)^{(1)}(t)$ in
$\Lambda_{(p_{i}-1)}(t)$.}
The following equivalent description of the subspaces $V_i(t)$ will
be very useful in the sequel:

\begin{lemma}
\label{complem}
A sequence of subspaces $\{\widetilde V_{\sigma_i}(t)\}_{i=1}^d$,
satisfying \eqref{esig02}, consists of the  canonical complements of
$\Bigl(\Lambda_{(p_{i})}\Bigr)^{(1)}(t)$ in $\Lambda_{(p_{i}-1)}(t)$
for any $1\leq i\leq d$ if and only if smooth (w.r.t. $t$) tuples of
vectors $E_{\sigma_i}(t)$, constituting bases of $\widetilde
V_{\sigma_i}(t)$, satisfy:
\begin{equation}
\label{normcond} \forall\, 1\leq j<i\leq d \text{ and } \forall\,
1\leq k\leq p_j-p_i+1:\quad
\omega\bigl(E_{\sigma_i}^{(p_i-1)}(t),E_{\sigma_j}^{(p_j-1+k)}(t)\bigr)=0
\end{equation}
or, equivalently, taking into account notations  in (\ref{esig6}),
\begin{equation}
\label{normcond2} 
 \forall\, 1\leq j<i\leq d \text{ and } \forall\,
1\leq k\leq p_j-p_i+1:\quad
\omega\bigl(E_{a_i}(t),E_{a_j}^{(k)}(t)\bigr)=0.
\end{equation}
\end{lemma}

The lemma can be easily proved by rewriting identity \eqref{cancomp}
in terms of bases $E_{\sigma_i}(t)$ and appropriate
differentiations.

Further, it turns out that on each canonical complement $V_i(t)$
one can define the canonical
quadratic form.
Indeed, given a vector $v\in V_i(t)$ take a smooth curve
$\varepsilon (t)$ in  $W$ such that
\begin{enumerate}
\item
$\varepsilon(\tau)=v$;
\item $\varepsilon(t)\in V_i(t)$ for any $t$ close to $\tau$.
\end{enumerate}
 Then by our constructions it is easy to see that for any
 $0\leq j\leq p_i-1$
\begin{subequations}
\label{vareps}
\begin{gather}
\label{vareps:a}\varepsilon^{(j)}(\tau)\in\Lambda_{(p_i-1-j)}(\tau),\\
\label{vareps:b}\varepsilon^{(j+1)}(\tau)\notin\Lambda_{(p_i-1-j)}(\tau),
\quad
 \text {if $v\neq 0$},\\
\label{vareps:c}\varepsilon^{(j+1)}(\tau)\in\Lambda_{(p_i-1-j)}(\tau),
\quad
 \text {if $v=0$}
 \end{gather}
 \end{subequations}
For this take a basis $E_{\sigma_i(t)}$ of $V_i(t)$, depending
smoothly on $t$,  expand our curve $\varepsilon(t)$ w.r.t. this
basis, and use the fact that for any $0\leq j\leq p_i-1$
\begin{equation}
\label{ee} \bigoplus_{s=0}^j {\rm span}\,
\{E_{\sigma_i}^{(s)}(t)\}\subset \Lambda_{(p_i-1-j)}(\tau),\quad
{\rm span}\, \{E_{\sigma_i}^{(j+1)}(t)\}\cap
\Lambda_{(p_i-1-j)}(\tau)=0,
\end{equation}
which is a direct consequence of Lemma \ref{filling}. From
(\ref{vareps:a}), (\ref{vareps:c}), the fact that  $\Lambda(t)$ is
the curve of Lagrangian subspaces, and the identity (\ref{subsup})
it follows that
\begin{equation}
\label{quadr} Q_{i,
\tau}(v)=\omega\bigl(\varepsilon^{(p_i)}(\tau),\varepsilon^{(p_i-1)}(\tau)\bigr)
\end{equation}
 is a well defined quadratic form on $V_i(\tau)$, which does not depend on the choice of the curve $\varepsilon(\tau)$ satisfying
conditions (1) and (2) above. The form  $Q_{i, \tau}(v)$ will be called
the \emph{canonical quadratic form on $V_i(\tau)$}. The quadratic
forms $Q_{i, \tau}(v)$ are nondegenerated for any $1\leq i\leq d$.
Indeed, if tuples $E_{\sigma_i(t)}$ constitute bases of $V_i(t)$ for
any $1\leq i\leq d$ and $Z_d(t)$ is as in \eqref{Z}, then from Lemma
\ref{complem} it follows that the matrix of the quadratic form
$\dot\Lambda(\tau)|_{Z_d(\tau)}$ in the basis
$\{E_{\sigma_k}^{(p_k-1)}(\tau)\}_{k=1}^d$ is block-diagonal and the
diagonal blocks are exactly the matrices of the forms $Q_{i,
\tau}(v)$ in the bases $E_{\sigma_i(t)}$. Then the nondegenericity
of the form $Q_{i, \tau}(v)$ follows from condition (G) and
\eqref{mon2}.
Moreover, if the curve $\Lambda(t)$ is monotonically nondecreasing,
then the forms $Q_{i, \tau}$ are positive definite. In this case the
Euclidean structure on $V_{\sigma_i}(\tau)$, corresponding to the
form $Q_{i, \tau}$ will be called the \emph{canonical Euclidean
structure on $V_i(\tau)$}.

From now on for simplicity of presentation we will assume that the
curve $\Lambda(t)$ is monotonically nondecreasing. All necessary
changes in the formulation of the results for nonmonotonic curves,
satisfying condition (G), will be indicated in section 4.
For any $1\leq i\leq d$, let ${\mathfrak B}_i$ be a fiber bundle over
the curve $\Lambda(t)$ such that the fiber of ${\mathfrak B}_i$ over
the point $\Lambda(t)$ consists of all
orthonormal bases of the space $V_{\sigma_i}(t)$ w.r.t. the
canonical Euclidean structure on $V_i(t)$.
Note that ${\mathfrak B}_i$ is the principle bundle with the structure group
$O(r_i)$.

\subsection{The canonical connections on the bundles $\mathfrak B_i$.} Now let us prove the following
\begin{prop}
\label{cancon} Each bundle ${\mathfrak B}_i$ is endowed with the
\emph{canonical principal connection} uniquely characterized by the
following condition: the section $E_{\sigma_i}(t)$ of ${\mathfrak
B}_i$ is horizontal w.r.t. this connection if and only if  ${\rm
span}\{E_{\sigma_i}^{(p_i)}(t)\}$ are isotropic subspaces of $W$ for
any $t$. Given any two horizontal sections $E_{\sigma_i}(t)$ and
$\widetilde E_{\sigma_i}(t)$ of $\mathfrak B_i$ there exists a
constant orthogonal matrix $U_i$
such that
\begin{equation}
\label{Uconst} \widetilde E_{\sigma_i}(t)=E_{\sigma_i}(t)U_i.
\end{equation}
\end{prop}

{\bf Proof}.
 As in the proof
of Lemma \ref{simpilem}, given two tuples of vectors $V_1=(v_{11},
\ldots, v_{1n_1})$ and $V_2=(v_{21}, \ldots, v_{2n_2})$ by
$\omega(V_1, V_2)$ we will mean the $n_1\times n_2$-matrix with the
$(i,j)$-entry equal to $\omega(v_{1i}, v_{2j})$. With this notation,
it is obvious that if $V_i={\rm span}\{ \widetilde E_{\sigma_i}\}$,
then the subspace ${\rm span}\{\widetilde
E_{\sigma_i}^{(p_i)}(t)\}$ is isotropic if and only if
\begin{equation}
\label{isot1} \omega\bigl(\widetilde E_{\sigma_i}^{(p_i)}(t),
\widetilde E_{\sigma_i}^{(p_i)}(t)\bigr)=0.
\end{equation}
Note also that from definition of the canonical Euclidean structure
it follows immediately that  for any section $ E_{\sigma_i}(t)$ of
the bundle $\mathfrak B_i$ the following identity holds
\begin{equation}
\label{isot15} \omega\bigl(E_{\sigma_i}^{(p_i)}(t),
E_{\sigma_i}^{(p_i-1)}(t)\bigr)=\rm {Id}.
\end{equation}
Take any two section $E_{\sigma_i}(t)$ and $\widetilde
E_{\sigma_i}(t)$ of the bundle $\mathfrak B_i$. Then there exists a
curve $U_i(t)$ of orthonormal matrices such that
$\widetilde E_{\sigma_i}(t)=E_{\sigma_i}(t)U_i(t)$.
Using relation
$\Lambda_{(1)}(t)=\bigl(\Lambda^{(1)}(t)\bigr)^\angle$ and formula
\eqref{isot15}, it is easy to get that
$$\omega\bigl(\widetilde E_{\sigma_i}^{(p_i)}(t), \widetilde
E_{\sigma_i}^{(p_i)}(t)\bigr)=U(t)^T\left(2p_i
U'(t)+\omega\bigl(E_{\sigma_i}^{(p_i)}(t),
E_{\sigma_i}^{(p_i)}(t)\bigr) U(t)\right).$$ So, relation
\eqref{isot1} holds if and only the matrix $U(t)$ satisfies the
following differential equation
\begin{equation}
\label{isot2} 2p_i U'(t)+\omega\bigl( E_{\sigma_i}^{(p_i)}(t),
E_{\sigma_i}^{(p_i)}(t)\bigr) U(t)=0.
\end{equation}
Note that the matrix $\omega\bigl( E_{\sigma_i}^{(p_i)}(t),
E_{\sigma_i}^{(p_i)}(t)\bigr)$ is antisymmetric. So, equation
\eqref{isot2} has solutions in $O(r_i)$, which are defined up to the
right translation there. This completes the proof of the
proposition.
$\Box$
\medskip

Now, if  for any $1\leq i\leq d$ we take a horizontal section
$E_{\sigma_i}(t)$ of the bundle ${\mathfrak B}_i$ and set , as
before, $E_{l^{j}(\sigma_i)}(t)=E_{\sigma_i}^{(j)}(t)$ for $0\leq j
\leq p_i-1$, then
from (\ref{Uconst}) it follows that for any superbox $a$  the
subspaces $V_a(t)={\rm span} \{E_a(t)\}$  do not depend on the
choice of a horizontal sections $E_{\sigma_i}(t)$.
Moreover, from this and Lemma \ref{filling} we get the \emph
{canonical splitting $\Lambda(t)=\bigoplus_{a\in \Delta} V_a(t)$ of
the subspaces $\Lambda(t)$}.

\subsection{The completion of horizontal sections to quasi-normal moving frames.}
In the sequel it will be more convenient to use the following
obviously equivalent description of quasi-normal mappings:

\begin{lemma}
\label{quasi2} A symmetric compatible mapping $R:\Delta\times
\Delta\mapsto \rm {Mat}$ is  quasi-normal  if and only if the
following four conditions hold: \vskip .1in
\begin{enumerate}

\item If $a$ and $b$ are two consecutive  superboxes in the same level of
$\Delta$, then the matrix $R(a,b)$ is antisymmetric; \vskip .1in

\item If both superboxes $a$ and $b$ are not special and
do not lie in the same or adjacent columns, then $R(a,b)=0$;\vskip
.1in

\item If both superboxes $a$ and $b$ are not special,
lie in the adjacent (but not the same) columns and one of the
superboxes is located from below and from  the left w.r.t. the
other, then $R(a,b)=0$;\vskip .1in

\item If a superbox $a$ is special, a superbox $b$ is not
special and $b$ is located from the left to $a$, but not in the
adjacent column, then $R(a,b)=0$.
\end{enumerate}
\end{lemma}

Further,  for all $1\leq i\leq d$, fix a horizontal section
$E_{\sigma_i}(t)$ of the bundle $\mathfrak B_i$ and complete it to
the moving basis $\{E_a(t)\}_{a\in\Delta}$ of $\Lambda(t)$ setting,
as before, $E_{l^{j}(\sigma_i)}(t)=E_{\sigma_i}^{(j)}(t)$ for $0\leq
j \leq p_i-1$. Also let
\begin{equation}
\label{fa} F_{a_i}(t)=E_{a_i}'(t).
\end{equation}
From the definition of the canonical Euclidean structure
it follows that $\omega(F_{a_i}(t), E_{a_i}(t))=Id$. From the normalization
conditions (\ref{normcond2}) with $k=1$ it follows that
$\omega(F_{a_i}(t), E_{a_j}(t))=0$ for any $i\neq j$. Further, by
definition of the horizontal section of the bundle $\mathfrak B_i$
one has $\omega(F_{a_i}(t), F_{a_i}(t))=0$. Finally, from the
normalization conditions (\ref{normcond2}) with $k=2$ it follows
that $\omega(F_{a_i}(t), F_{a_j}(t))=0$ for $i\neq j$ as well.
Combining all these identities with the fact that the subspaces
$\Lambda(t)$ are Lagrangian and the relation
$\Lambda_{(1)}(t)=\bigl(\Lambda^{(1)}(t)\bigr)^\angle$, we get that
the tuple $\bigl(\{E_a\}_{a\in \Delta}, \{F_b(t)\}_{b\in{\mathcal
F_1}})$, where, as before, $\mathcal F_1$ denotes the first column
of $\Delta$, does not contradict the relations for a Darboux frame.
Besides, by our constructions it satisfies first two equations of
(\ref{structeq}). In this subsection we prove the following

\begin{prop}
\label{quasilem} The tuple $\bigl(\{E_a\}_{a\in \Delta},
\{F_b(t)\}_{b\in{\mathcal F_1}})$ can be uniquely completed to a
quasi-normal moving frame of the curve $\Lambda(t)$.
\end{prop}

{\bf Proof}.
 Take a tuple $\{
F_b(t)\}_{b\in \Delta\backslash\mathcal F_1}$, which completes the
tuple $\bigl(\{E_a\}_{a\in \Delta}, \{F_b(t)\}_{b\in{\mathcal
F_1}})$ to a moving Darboux's frame in $W$. Then from the definition
of Darboux's frame
and the first two equations of (\ref{structeq}) it follows that this
moving Darboux frame have the structural equation (\ref{structeq})
for some symmetric mappings $
R_t:\Delta\times\Delta\mapsto {\rm Mat}$ compatible with the Young
diagram $D$.
As before, denote by $\mathcal F_j$ the $j$th column of $\Delta$,
$1\leq j\leq p_1$. Our proposition will follow from the following

\begin{stat} For any $1\leq k\leq p_1$ there exists a unique tuple of
columns of vectors $$\{F_b(t): b\in \bigcup_{j=1}^{k}\mathcal
F_j\}$$ such that the tuple $\bigl(\{E_a\}_{a\in \Delta},\{F_b(t):
b\in \bigcup\limits_{j=1}^{k}\mathcal F_j\})$ can be completed to a
moving Darboux frame $\bigl(\{E_a\}_{a\in \Delta},
\{F_b(t)\}_{b\in\Delta})$ such that if the mapping
$R_t:\Delta\times\Delta\mapsto{\rm Mat}$ appears in the structural
equation  \eqref{structeq} for this moving frame, then the mapping
$R_t$ satisfies conditions (1)-(4) of Lemma \ref{quasi2} for any
pair $(a,b)$ with at least one superbox belonging to the first
$(k-1)$ columns of $\Delta$.
\end{stat}

Indeed, our proposition is just Statement 1 in the case $k=p_1$ (the
only pair of superboxes, which is not covered by Statement 1, is
$(\sigma_1,\sigma_1)$, where, as before, $\sigma_1$ is the special
(the last) superbox of the first level, but this pair does not
satisfy any of conditions (1)-(4) of Lemma \ref{quasi2}).

We will prove Statement 1 by induction w.r.t. $k$. For $k=1$ there is
nothing to prove, because the tuple $\{F_c\}_{c\in\mathcal F_1}$ is
uniquely determined by the second line of (\ref{structeq}) (which
together with the first line of (\ref{structeq}) is equivalent to
(\ref{fa})), while the Statement 1 for $k=1$ does not impose any
conditions on the symmetric compatible mapping $R_t$, appearing in
(\ref{structeq}).

Now suppose that Statement 1 is proved for some $k=\bar k$, where
$1\leq \bar k\leq p_1-1$, and prove it for $k=\bar k+1$. Let
$\{F_b(t): b\in \bigcup\limits_{j=1}^{\bar k}\mathcal F_j\}$ be the
tuple, satisfying Statement 1 for $k=\bar k$.
Take a tuple $\{
F_b(t):b\in \Delta\backslash\bigcup\limits_{j=1}^{\bar k}\mathcal
F_j\}$, which completes the tuple $\bigl(\{E_a\}_{a\in \Delta},
\{F_b(t): b\in \bigcup\limits_{j=1}^{\bar k}\mathcal F_j\})$ to a
moving Darboux's frame in $W$ and assume that
$R_t:\Delta\times\Delta\mapsto{\rm Mat}$ is the mapping, appearing
in the structural equation for this frame. If $\{\widehat F_b(t):
b\in \Delta\backslash\bigcup\limits_{j=1}^{\bar k}\mathcal F_j\}$ is
another tuple, completing the tuple $\bigl(\{E_a\}_{a\in \Delta},
\{F_b(t)\}_{b\in{\mathcal F_1}})$ to a moving Darboux's frame in
$W$,
then there exists a symmetric mapping
$\Gamma_t:(\Delta\backslash\bigcup\limits_{j=1}^{\bar k}\mathcal
F_j)\times (\Delta\backslash\bigcup\limits_{j=1}^{\bar k}\mathcal
F_j)\mapsto {\rm Mat}$, compatible with the diagram, obtained from
$D$ by erasing the first $\bar k$ column, such that
\begin{equation}
\label{transDb} \forall a\in
\Delta\backslash\bigcup\limits_{j=1}^{\bar k}\mathcal F_j\quad
\widehat F_a(t)=
F_a(t)+ \sum_{b\in \Delta\backslash\bigcup\limits_{j=1}^{\bar
k}\mathcal F_j}E_b(t)\Gamma_t(a,b).
\end{equation}
Suppose that $\widehat R_t:\Delta\times\Delta\mapsto {\rm Mat}$ is
the symmetric mapping compatible with the Young diagram $D$ such
that similarly to last two equations of
(\ref{structeq}) one has
\begin{equation}
\label{structeq2}
\begin{cases}
F_a'(t)=\sum\limits_{b\in\Delta}\widehat E_b R_t(a,b)-\widehat
F_{r(a)}&\text{if $a\in \bigcup\limits_{j=1}^{\bar k}\mathcal F_j$}\\
 \widehat F_a'(t)=\sum\limits_{b\in\Delta}\widehat
E_b R_t(a,b)-\widehat F_{r(a)}& \text {if\,\, $a\in
\Delta\backslash(\bigcup\limits_{j=1}^{\bar k}\mathcal F_j\cup\mathcal S)$}\\
\widehat F_a'(t)=\sum\limits_{b\in\Delta}\widehat E_b R_t(a,b)&
\text {if\,\, $a\in \mathcal S$},
\end{cases}
\end{equation}
(note that from the first line of (\ref{structeq2}), one has
$\widehat R_t(a,b)=R_t(a,b)$, if at least one of the superboxes
$(a,b)$ belongs to the first $\bar k$ columns of $\Delta$).
Let us extend the mappings
$\Gamma_t:(\Delta\backslash\bigcup\limits_{j=1}^{\bar k}\mathcal
F_j)\times (\Delta\backslash\bigcup\limits_{j=1}^{\bar k}\mathcal
F_j)\mapsto {\rm Mat}$ to the symmetric mapping, denoted by the same
letter $\Gamma_t$, from $\Delta\times\Delta$ to ${\rm Mat}$
compatible with the diagram $D$, by setting
\begin{equation}
\label{bm} \Gamma_t(a,b)(t)=\Gamma_t(b,a)^T=0,\quad \forall\, b\in
\bigcup\limits_{j=1}^{\bar k}\mathcal F_j, a\in\Delta.
\end{equation}
Then, substituting (\ref{transDb}) into two last lines of
(\ref{structeq2}) and using (\ref{structeq}), one can easily obtain
\begin{equation}
\label{transxi}
\widehat R_t(a,b)=R_t(a,b)+\frac{d}{dt}\Gamma_t(a,b)+
\Gamma_t\bigl(a, r(b)\bigr)+\Gamma_t
\bigl(r(a),b\bigr),\end{equation} where the term $\Gamma_t(a, r(b))$
is omitted, if $b$ is special, and the term $\Gamma_t (r(a),b)$ is
omitted, if $a$ is special.
 Using transformation rule
(\ref{transxi}), we will prove the following

\begin{stat}
There exists the unique choice of matrices $\Gamma_t\bigl(\tilde a,
\tilde b\bigr)$ with at least one of the superboxes belonging to the
$(\bar k+1)$th column of $\Delta$ and the other one lying from the
right to the $\bar k$th column of $\Delta$ such that the matrix
$\widehat R_t(a,b)$ satisfies all conditions (1)-(4) of Lemma
\ref{quasi2} for any pairs $(a,b)$ with at least one of the
superboxes belonging to the $\bar k$th column of $\Delta$ and the
other one lies from the right to the $(\bar k-1)$th column of
$\Delta$
\end{stat}
It is clear that Statement 2, relation (\ref{transDb}), and the
induction hypothesis will imply Statement 1 for $k=\bar k+1$.
Let us prove statement 2.  Suppose that $a\in \mathcal F_{\bar k}$.
Then from (\ref{bm}) it follows that $\frac{d}{dt}\Gamma_t(a,b)=0$
and $\Gamma_t\bigl(a, r(b)\bigr)=0$. So, relations (\ref{transxi})
in this case have a form
\begin{equation}
\label{transxi1}
\widehat R_t(a,b)=
R_t(a,b)+\Gamma_t\bigl(r(a), b\bigr),\end{equation} where the term
$\Gamma_t\bigl(r(a), b\bigr)$ is omitted, if $a$ is special
(obviously it happens, when the level of $a$ consists of only one
superbox). Therefore, according to (\ref{transxi1}), if $a$ is special or
$b\in \bigcup\limits_{j=1}^{\bar k}\mathcal F_j$ we have $\widehat
R_t(a,b)= R_t(a,b)$, i.e. the matrix $R_t(a,b)$ is already
independent of the choice of the complement of $\bigl(\{E_{\tilde
a}\}_{\tilde a\in \Delta}, \{F_{\tilde b}(t):\tilde b\in
\bigcup\limits_{j=1}^{\bar k}\mathcal F_j\}\bigr)$ to a moving
Darboux frame.

Now assume that $a$ is not special and $b\notin
\bigcup\limits_{j=1}^{\bar k}\mathcal F_j$. Then there are the
following  three cases:

{\bf a)} $b\notin \bigcup\limits_{j=1}^{\bar k+1}\mathcal F_j$, i.e.
$b$ is not in the first  $\bar k+1$ columns of $\Delta$.
Then the matrix $\Gamma_t\bigl(r(a), b\bigr)$ appears only ones in
all relations, \begin{equation} \label{transxi1t}
\widehat R_t(\tilde a,\tilde b)=
R_t(\tilde a,\tilde b)+\Gamma_t\bigl(r(\tilde a), \tilde
b\bigr),\end{equation} where $\tilde a$ runs over the whole $\bar
k$th column $\mathcal F_{\bar k}$ of $\Delta$. Putting
\begin{equation}
\label{transxi2} \Gamma_t\bigl(r(a), b\bigr)=-
R_t(a,b),
\end{equation}
we get $\widehat R_t(a,b)=0$ for any $a\in{\mathcal F_{\bar
k}}$, which  corresponds to conditions (2) and (4) of Lemma \ref{quasi2}, if
$b$ is not from the left to $a$.
Obviously, the choice of $\Gamma_t\bigl(r(a), b\bigr)$ as in
(\ref{transxi2}) is the unique one with these properties.

{\bf b)} $b\in \mathcal F_{\bar k+1}$, but $b\neq r(a)$, i.e. $b$
lies in the $(\bar k+1)$th column of $\Delta$, but it is not in the
same row with $a$. Let $a_1=l(b)$. Then from the symmetricity of the
mapping $\Gamma_t$ (i.e. the relation
$\Gamma_t(a,b)=\bigl(\Gamma_t(a,b)\bigr)^T$) it follows that the
matrix $\Gamma_t\bigl(r(a_1), r(a)\bigr)$ appears twice in all
relations (\ref{transxi1t}), where $\tilde a$ runs over the $\bar
k$th column $\mathcal F_{\bar k}$ of $\Delta$ and $\tilde b$ runs
over the $(\bar k+1)$th column $\mathcal F_{\bar k+1}$ of $\Delta$.
Namely, substituting $(\tilde a, \tilde b)=\bigl(r(a), a_1\bigr)$
into (\ref{transxi1t}) and using the symmetricity of the mapping
$\Gamma_t$  we will get the following relation in addition to
(\ref{transxi1}) (with $b=r(a_1)$):
\begin{equation}
\label{transxi3}
\widehat R_t\bigl(a_1, r(a)\bigr)= R_t\bigl(a_1,
r(a)\bigr)+\Gamma_t\bigl(r(a), r(a_1)\bigr)^{T}.\end{equation}
Hence, from symmetricity again we have
$$\widehat R_t\bigl(a,r(a_1)\bigr)-\widehat R_t\bigl(r(a), a_1\bigr)=
R_t\bigl(a, r(a_1)\bigr)-
R_t\bigl(r(a),a_1\bigr),$$ i.e. the matrix  $
R_t\bigl(a,r(a_1)\bigr)-
R_t\bigl(r(a), a_1\bigr)$ does not depend on the choice of the
complement of $\bigl(\{E_{\tilde a}\}_{\tilde a\in \Delta},
\{F_{\tilde b}(t):\tilde b\in \bigcup\limits_{j=1}^{\bar k}\mathcal
F_j\}\bigr)$ to a moving Darboux frame. Besides, for any pair of
superboxes $(a,a_1)$, $a\neq a_1$ in the $\bar k$th column $\mathcal
F_{\bar k}$ by an appropriate choice of $\Gamma_t\bigl(r(a),
r(a_1)\bigr)$ we cannot "kill" both matrices $
R_t\bigl(r(a),a_1\bigr)$ and
$
R_t\bigl(a,r(a_1)\bigr)$, but only one of them.
We choose
the following normalization: $\widehat R(a, r(a_1))=0$, if $a_1$ is
higher than $a$. We can do it by putting $\Gamma_t\bigl(
r(a),r(a_1))\bigr)=-
R_t\bigl(a,r(a_1)\bigr)$. This normalization corresponds to conditions (3) of Lemma
\ref{quasi2}.
Obviously, such choice of $\Gamma_t\bigl(r(a),
r(a_1)\bigr)$ is the unique one with these properties.

{\bf c)} $b=r(a)$. Then the matrix $\Gamma_t\bigl(r(a), r(a)\bigr)$
appears only once in all relations \eqref{transxi1t} where $\tilde
a$ runs over the whole $\bar k$th column $\mathcal F_{\bar k}$ of
$\Delta$, namely
\begin{equation} \label{transxi1s}
\widehat R_t(a,r(a))=R_t(a, r(a))+\Gamma_t\bigl(r(a), r(
a)\bigr).\end{equation}
 On
the other hand, by our assumptions $\Gamma_t\bigl(r(a), r(a)\bigr)$
should be symmetric. Therefore,  using (\ref{transxi1s}), we cannot
"kill" the whole matrix $
R_t(a,r(a))$, but only its symmetric
part (by putting $\Gamma_t\bigl(r(a), r( a)\bigr)=-\frac
{1}{2}\bigl(
R_t(a, r(a))+
R_t(a, r(a))^T\bigr)$). It corresponds to conditions (1) of Lemma
\ref{quasi2} with $a\in{\mathcal F_{\bar k}}$. Obviously, such
choice of $\Gamma_t\bigl(r(a), r(a)\bigr)$ is the unique one with
these properties. In this way we have found uniquely all matrices
$\Gamma_t(\tilde a,\tilde b)$ with $\tilde a \in \mathcal F_{\bar k
+1}$, $b\notin \bigcup\limits_{j=1}^{\bar k}\mathcal F_j$ such that
the matrix $\widehat R_t(a,b)$ satisfies all conditions (1)-(4) of
Lemma \ref{quasi2} for any pairs $(a,b)$, where $a \in \mathcal
F_{\bar k}$, $b\notin \bigcup\limits_{j=1}^{\bar k-1}\mathcal F_j$.
Taking $\Gamma_t(\tilde b,\tilde a)=\Gamma_t(\tilde a, \tilde b)^T$,
we will have the same properties for $\widehat R_t(b,a)$ with $a$
and $b$ as in the previous sentence. This completes the proof of
Statement 2, therefore also the proof of the Statement 1 for $k=\bar
k+1$, and then by induction the proof of Proposition \ref{quasilem}.
$\Box$
\medskip

\subsection{Normality of the obtained quasi-normal moving frames}
In the present subsection we will show that the quasi-normal moving
frame, constructed in the previous subsection, is in fact a normal
moving frame. Note that in the previous subsection we did not use at
all the normalization conditions (\ref{normcond2}) with $k\geq 3$.
As before, we denote by $d$ the number of levels in the diagram
$\Delta$,  by $p_i$ the number of superboxes in the $i$th level, and
by $a_i$ the first superbox in the $i$th level. The normality of the
constructed quasinormal frame will obviously follow from the
following

\begin{prop}
\label{extra0} A quasi-normal moving frame $(\{E_a(t)\}_{a\in
\Delta}, \{F_a(t)\}_{a\in \Delta})$ is normal if and only if
conditions (\ref{normcond2}) hold for any $1\leq j<i\leq d$ and
$3\leq k\leq p_j-p_i+1$.
\end{prop}

Proposition \ref{extra0} will follow by induction from the following

\begin{stat}
\label{stat6} Fix $s\in \mathbb N$ and let $R_t:\Delta \times
\Delta\to {\rm Mat}$ be a quasi-normal mapping, satisfying the
following condition: for
any $i$ and  $j$, $1\leq j<i\leq d$, the matrix $R_t(a,b)\equiv 0$
for all first $\min\{s-1,p_j-p_i-1\}$ pairs $(a,b)$ in the tuple
\eqref{chain}. Then for any $i$ and $j$, $1\leq j<i\leq d$, such
that $1\leq s\leq p_j-p_i$, the $s$th pair $(\bar a_i^s, \bar
a_j^s)$ of the tuple \eqref{chain} satisfies
\begin{equation}
\label{spair} R_t(\bar a_i^s, \bar a_j^s)=\pm
\omega\bigl(E_{a_j}^{(s+2)}(t), E_{a_i}(t)\bigr).
\end{equation}
\end{stat}

Before proving Statement \ref{stat6}, let us introduce some
notations. As in the proof of Lemma \ref{simpilem}, let
$\overline\Delta$ be the diagram obtained from $\Delta$ by the
reflection w.r.t. its left edge. In the sequel we will work with the
diagram $\Delta\cup\overline\Delta$.
The boxes of this
diagram will be also called superboxes. Similar to above, we will
denote by $l$ and $r$ the left and the right shifts on the diagram
$\Delta\cup\overline\Delta$.
\begin{defin}
\label{admpath}
 A (finite) sequence $\eta=\{b_0,\ldots, b_n\}$ of superboxes of the diagram $\Delta\cup\overline\Delta$
is called an admissible path in this diagram, if the following two
conditions hold:
\begin{enumerate}
\item If $b_i\in \Delta$ then $b_{i+1}\in\{b_i,l(b_i)\}$;
\item If $b_i\in \overline \Delta$ then
$b_{i+1}\in\{b_i,l(b_i)\}\cup\Delta$
\end{enumerate}
{\rm (see an example on Figure 1)}. The superboxes from the admissible
path $\eta$ will be called the \emph{vertices} of the path. We will
distinguish three types of vertices: the vertex $b_m$, $0\leq m<n$,
will be called \emph{walking}, if $b_{m+1}=l(b_m)$, it will be
called \emph {sleeping}, if $b_{m+1}=b_m$, and it will be called
\emph{jumping}, if $b_m\in \overline\Delta$ and $b_{m+1}\in\Delta$.
\end{defin}
\setlength{\unitlength}{0.65mm}
\begin{picture}
(130,50)(-40,15)
\put(20,56){\line(1,0){30}} \put(60,56){\line(1,0){20}}
\put(90,56){\line(1,0){30}} \put(20,56){\line(0,-1){7}}
\put(20,49){\line(1,0){30}}
\put(60,49){\line(1,0){20}} \put(90,49){\line(1,0){30}}

\put(30,56){\line(0,-1){7}} \put(40,42){\line(1,0){10}}
\put(60,42){\line(1,0){20}}
\put(90,42){\line(1,0){10}}\put(90,42){\line(0,1){10}}

\put(40,56){\line(0,-1){14}} \put(55,52) {\makebox(0,0){\ldots}}
\put(55,45) {\makebox(0,0){\ldots}}
\put(55,40) {\makebox(0,0){$\ddots$}} \put(60,35){\line(1,0){20}}
\put(60,28){\line(1,0){20}} \put(80,56){\line(0,-1){28}}
\put(90,56){\line(0,-1){14}} \put(100,56){\line(0,-1){14}}
\put(110,56){\line(0,-1){7}}
\put(120,56){\line(0,-1){7}}\put(50,56){\line(0,-1){14}}
\put(70,56){\line(0,-1){28}} \put(70,56){\line(0,1){3}}
\put(70,28){\line(0,-1){3}} \put(85,52) {\makebox(0,0){\ldots}}
\put(85,45) {\makebox(0,0){\ldots}}
 \put(83 ,37) {\makebox(0,0){.}}
 \put(85,39) {\makebox(0,0){.}}
\put(87 ,41) {\makebox(0,0){.}} \put(60,56){\line(0,-1){28}}

\put(70,18){\makebox(0,0){Figure 1.}}
\put(60,56){\makebox(5,10){$\overline \Delta$}}
\put(70,56){\makebox(15,9){$\Delta$}}

\put(105,53){\circle*{2}} \put(105,58){\makebox(0,0){$b_0$}}
\put(105,53){\line(-1,0){10}} \put(95,53){\circle*{2}}
\put(95,58){\makebox(0,0){$b_1$}} \put(95,53){\line(-1,0){7}}
\put(81.5,53){\line(-1,0){6.5}}\put(75,53){\circle*{2}}
\put(75,53){\line(-1,0){10}}\put(65,53){\circle*{2}}
\put(65,53){\line(-1,0){7}} \put(51.5,53){\line(-1,0){6.5}}
\put(45,53){\circle*{2}} \put(45,53){\line(6,-1){51}}
\put(95,44.5){\circle*{2}} \put(95,44.5){\line(-1,0){7}}
\put(81.5,44.5){\line(-1,0){6.5}}\put(75,44.5){\circle*{2}}
\put(75,44.5){\line(-1,0){10}}\put(65,44.5){\circle*{2}}
\put(65,44.5) {\line(2,-3){9.3}} \put(75,30.5){\circle*{2}}
\put(75,37.2){\makebox(0,0){$b_n$}}
%
\end{picture}

Further, given any superbox $x$ of $\Delta\cup\overline\Delta$ we
will denote by $\bar x$ the superbox , obtained from $x$ by the
reflection of $x$ w.r.t. the left edge of the diagram $\Delta$. We
also assume that the size of the superbox $x\in \overline \Delta$ is
equal to the size of superbox $\bar x$.

From the definition of Darboux frame it follows that the quantity
$-\omega( E_{a_i}, E_{a_j}^{(s+2)})$, we are interested in, is equal
to the coefficient near $F_{a_i}$ of the expansion of
$E_{a_j}^{(s+2)}$ into linear combination w.r.t. the frame
$(\{E_a(t)\}_{a\in \Delta}, \{F_a(t)\}_{a\in \Delta})$, satisfying
the structural equation (\ref{structeq}). Admissible pathes in the
diagram $\Delta\cup\overline\Delta$ help to describe the
coefficients of such expansions.
For this to any admissible path $\eta=\{b_0,\ldots, b_n\}$ we will
assign a curve of $\text{size}(b_n)\times\text{size}(b_0)$-matrices
$P_{\eta}(\cdot)$. The curve of matrices $P_\eta(\cdot)$ can be
defined by the recursive formulas on the number of vertices in
$\eta$. If $\eta$ consists of only one vertex, $\eta=\{b_0\}$, we
set $P_\eta(t)$ to be the identity matrix for any $t$. Further for
the path $\eta=\{b_0,\ldots, b_{n-1}, b_n\}$ the curve of matrices
$P_\eta(\cdot)$ is obtained  from the curve of matrices
$P_{\{b_0,\ldots, b_{n-1}\}}$ by the following recursive formula:

\begin{equation}
\label{derpath} P_{\{b_0,\ldots, b_{n-1}, b_n\}}(t)=\begin{cases}
P_{\{b_0,\ldots,
b_{n-1}\}}(t)&\text{if}\,\, b_n=l(b_{n-1}),\,\, b_{n-1}\in\Delta,\\
-P_{\{b_0,\ldots, b_{n-1}\}}(t)&\text{if}\,\, b_n=l(b_{n-1}),\,\,
b_{n-1}\in\overline\Delta,\\
P_{\{b_0,\ldots, b_{n-1}\}}^\prime(t)&\text{if}\,\, b_n=b_{n-1}, \\
R_t(\bar b_{n-1}, b_n) P_{\{b_0,\ldots,b_{n-1}\}}(t)&\text{if}\,\,
b_{n-1}\in\overline \Delta , \,\, b_n\in\Delta
\end{cases}
\end{equation}

Given $\{a,b\}\subset \Delta\cup\overline\Delta$
 and $n\in\mathbb N\cup\{0\}$ denote by $\Omega(a,b,n)$  the set of all admissible pathes in the
diagram $\Delta\cup\overline\Delta$,
 starting at $a$,
ending at $b$, and consisting of $n+1$ vertices. Then from
structural equation \eqref{structeq}, definition (\ref{derpath}) of
matrices $P_\eta$, and elementary rules of differentiations it
follows that \begin{equation} \label{admfix} \omega( E_{a_i},
E_{a_j}^{(s+2)})=-\sum_{\eta\in \Omega(a_j,\bar a_i,s+2)}P_\eta
\end{equation}

\begin{remark}
\label{jumping} \rm {It is clear from the last line of the recursive
formula (\ref{derpath}) that if $P_\eta(t)\neq 0$, then $R_t(\bar
b_m, b_{m+1})\neq 0$ for any jumping vertex $b_m$ of $\eta$.} $\Box$
\end{remark}

Further , it is convenient to enumerate the columns of the diagram
$\Delta\cup\overline\Delta$ by integers in the following way: to the
$j$th column (from the left) of $\Delta$  we assign the same number
$j$ while to the $j$th column from the right of $\overline{\Delta}$
we assign the number $1-j$. Given a superbox $a\in \Delta\cup
\bar{\Delta}$, denote by $c(a)$ the number of the column, according
to the rule described in the previous sentence. The following simple
lemma will be useful in the sequel

\begin{lemma}
\label{triv} Suppose that $R_t:\Delta\times \Delta\mapsto {\rm Mat}$
is a quasi-normal mapping and $R_t(a,b)\neq 0$, where superboxes $a$
and $b$ lie in the $j$th and $i$th level of $\Delta$ respectively
($j<i$). Then the pair $(a,b)$ is $\bigl(c(b)-c(\bar a)\bigr)$th
pair in the tuple (\ref{chain}).
\end{lemma}

Indeed, 
by Definition \ref{qnormmap} the nonzero matrix $R_t(a,b)$ must
correspond to a pair from the appropriate tuple of the form
(\ref{chain}). The second sentence of the lemma is obvious.

 {\bf Proof of Statement \ref{stat6}.} Fix some admissible path
$\eta=\{b_0,\ldots, b_{s+2}\}$ from $\Omega(a_j, \bar a_i, s+2)$ (by
definition, $b_0=a_j$ and  $b_{s+2}=\bar a_i$). Let us denote by $k$
the number of jumping vertices in $\eta$. Further, let
$b_{m_1},\ldots, b_{m_k}$ be all jumping vertices of $\eta$, where
$m_1<m_2<\ldots<m_k$. Set also $m_0=-1$, $m_{k+1}=s+2$. It is
evident that for any $1\leq u\leq k+1$ the number of superboxes
between $b_{m_{u-1}+1}$ and $b_{m_{u}}$ (including $b_{m_{u-1}+1}$
but not $b_{m_{u}}$) is equal to
$c\bigl(b_{m_{u-1}+1}\bigr)-c\bigl(b_{m_{u}}\bigr)$. Therefore the
fact that all superboxes $b_u$ with $0\leq u< s+1$ are either
walking or sleeping or jumping can be expressed as follows

\begin{equation}\label{times}
 \sum_{u=1}^{k+1}\bigl(c(b_{m_{u-1}+1})-c (b_{m_{u}})\bigr)+\#\{\text{sleeping vertices of }\eta\}+k=s+2.
\end{equation}

\begin{lemma}
\label{numjump} Under assumptions of Statement \ref{stat6} if
$P_\eta\neq 0$ for a path $\eta\in \Omega(a_j, \bar a_i, s+2)$
($j<i$) with
$p_j-p_i\geq s$,
then  there is only one jumping vertex and there are no sleeping
vertices in $\eta$.
\end{lemma}

{\bf Proof.}
 Since any path $\eta\in \Omega(a_j, \bar a_i, s+2)$ has
to contain at least one jumping vertex (in order to jump somehow
from $j$th to $i$th level) the lemma is actually equivalent to the
fact that
\begin{equation}
\label{slj1} \#\{\text{sleeping vertices of }\eta\}+k=1
\end{equation}

Assume the converse, i.e.
\begin{equation}
\label{slj2} \#\{\text{sleeping vertices of }\eta\}+k\geq 2.
\end{equation}

Given a superbox $x\in \Delta$, denote by $p(x)$ the number of
superboxes in the level of $x$. Assume that the superboxes
$b_{m_{u}}$ and $b_{m_{u}+1}$ lie in different levels. By Remark
\ref{jumping}, $R_t(\bar b_{m_{u}},b_{m_{u}+1})\neq 0$. Therefore,
according to Lemma \ref{triv} either $(\bar b_{m_{u}},
b_{m_{u}+1})$ or $( b_{m_{u}+1}, \bar b_{m_{u}})$ is the
$\bigl(c(b_{m_{u}+1})-c(b_{m_{u}})\bigr)$th pair in the tuple
(\ref{chain}). Combining this with Remark \ref{jumping} and
assumptions of Statement \ref{stat6}, one can obtain that if the
superboxes $b_{m_{u}}$ and $b_{m_{u}+1}$ lie in different levels,
then
\begin{equation}
\label{jump1}
c(b_{m_{u}+1})-c(b_{m_{u}})>\min\{s-1,|p(b_{m_{u}+1})-p(\bar
b_{m_{u}})|-1\}.
\end{equation}

Further, since  $c(b_0)=1$ and $c(b_{s+2})=0$ (recall that
$b_0=a_j$, $b_{s+2}=a_i$, and $m_{k+1}=s+2$), we have
\begin{equation}
\label{jump12} \sum_{u=1}^{k+1}\bigl(c(b_{m_{u-1}+1})-c
(b_{m_{u}})\bigr)= \sum_{u=1}^{k}\bigl(c(b_{m_{u}+1})-c
(b_{m_{u}})\bigr)+1.
\end{equation}
 Substituting the last identity into
(\ref{times}) and using assumption (\ref{slj2}) we obtain
\begin{equation}
\label{jump15}
\sum_{u=1}^{k}\bigl(c(b_{m_{u}+1})-c (b_{m_{u}})\bigr)
\leq s-1.
\end{equation}
 Since all terms in the sum in the lefthand side of the previous
inequality are positive, we have $c(b_{m_{u}+1})-c (b_{m_{u}})\leq
s-1$ for any $1\leq u\leq k$. Combining the last inequality with
(\ref{jump1}) we obtain that if the superboxes $b_{m_{u}}$ and
$b_{m_{u}+1}$ lie in different levels, then

\begin{equation}
\label{jump2}
c(b_{m_{u}+1})-c(b_{m_{u}})\geq|p(b_{m_{u}+1})-p(\bar b_{m_{u}})|.
\end{equation}
Besides, if the superboxes $b_{m_{u}}$ and $b_{m_{u}+1}$ lie in the
same level, then the inequality (\ref{jump2}) holds automatically.

On the other hand, by our constructions the superboxes $b_{m_{u}+1}$
and $\bar b_{m_{u+1}}$ lie in the same level of $\Delta$. This fact
together with inequalities (\ref{jump2}) and (\ref{jump15}) implies
that
$$p_j-p_i\leq \sum_{i=1}^{k}|p(b_{m_{u}+1})-p(\bar b_{m_{u}})| \leq
\sum_{i=1}^{k}c(b_{m_{u}+1})-c(b_{m_{u}})\leq s-1,$$ which
contradicts the assumption $p_j-p_i\geq s$ of Lemma \ref{numjump}.
The proof of the lemma is completed. $\Box$
\medskip

Now, if $\eta$ has only one jumping vertex and no sleeping vertices,
then from (\ref{times}) and (\ref{jump12}) it follows that
$c(b_{{m_1}+1})-c(b_{m_1})=s$. Besides, in this case the superbox
$b_{m_1}$ lies in the $j$th level and the superbox $b_{m_1+1}$ lies
in the $i$th level. But then from Remark \ref{jumping} and Lemma
\ref{triv} it follows that if $P_\eta\neq 0$ then the pair $(\bar
b_{m_1}, b_{{m_1}+1})$ is exactly the $s$th pair of the tuple
(\ref{chain}), which together with (\ref{derpath}) and
(\ref{admfix}) implies (\ref{spair}). The proof of Statement
\ref{stat6} is completed. $\Box$
\medskip

As we have already menstioned, Proposition \ref{extra0} follows
immediately from Statement \ref{stat6} by induction w.r.t. $s$,
starting with $s=1$ (for which the assumptions of Statement
\ref{stat6} hold automatically).
\medskip

\subsection{Final steps of the proof of Theorem \ref{maintheor}}
The "if" part of Proposition \ref{extra0} implies that the tuple
$(\{E_a(t)\}_{a\in \Delta}, \{F_a(t)\}_{a\in \Delta})$ constructed
in the subsection 3.5 is a normal moving frame of the curve
$\Lambda(\cdot)$. Moreover, by the constructions of subsection 3.3
the space $V_i(t)={\rm span}\{ E_{\sigma_i}(t)\}$ is the canonical
complement of $\Bigl(\Lambda_{(p_{i})}\Bigr)^{(1)}(t)$ in
$\Lambda_{(p_{i}-1)}(t)$ for any $1\leq i\leq d$, where $\sigma_i$
is the special superbox of the $i$th level, and by constructions of
subsection 3.4 the curves $E_{\sigma_i}(t)$ are horizontal sections
of the bundle $\mathfrak B_i$, defined in subsection 3.3.

Now suppose that $(\{\widetilde E_a(t)\}_{a\in \Delta}, \{\widetilde
F_a(t)\}_{a\in \Delta})$ is another normal moving frame of the curve
$\Lambda(\cdot)$. From the second line of the structural equation
(\ref{structeq}) (where all $E_a(t)$ and $F_a(t)$ are replaced by
$\widetilde E_a(t)$ and $\widetilde F_a(t)$) and the definition of
Darboux frame it follows that conditions (\ref{normcond2}) (again
with all $E_a(t)$ replaced by $\widetilde E_a(t)$) hold for any
$1\leq j<i\leq d$ and $k=1,2$. Indeed, $\omega\bigl(\widetilde
E_{a_i}(t), \widetilde E_{a_j}'(t)\bigr)=\omega\bigl(\widetilde
E_{a_i}(t), \widetilde F_{a_j}(t)\bigr)=0$ and
$\omega\bigl(\widetilde E_{a_i}(t), \widetilde
E_{a_j}''(t)\bigr)=-\omega\bigl(\widetilde E_{a_i}'(t), \widetilde
E_{a_j}'(t)\bigr)= -\omega\bigl(\widetilde F_{a_i}(t), \widetilde
F_{a_j}(t)\bigr)=0$. Further, by Proposition \ref{extra0}, from the
normality of the frame $(\{\widetilde E_a(t)\}_{a\in \Delta},
\{\widetilde F_a(t)\}_{a\in \Delta})$ it follows that conditions
(\ref{normcond2}) (again with all $E_a(t)$ replaced by $\widetilde
E_a(t)$) hold for any $1\leq j<i\leq d$ and $3\leq k\leq p_j-p_i+1$.
Therefore, Lemma \ref{complem}
implies that ${\rm span}\{ \widetilde E_{\sigma_i}(t)\}={\rm
span}\{ E_{\sigma_i}(t)\}=V_i(t)$. Besides, from the second line of
the structural equation (\ref{structeq}) (where again all $E_a(t)$
and $F_a(t)$ are replaced by $\widetilde E_a(t)$ and $\widetilde
F_a(t)$) and Proposition \ref{cancon}
it follows that the curves $\widetilde E_{\sigma_i}$ are horizontal
sections of the bundle $\mathfrak B_i$, which together with
(\ref{Uconst})
implies relations (\ref{u}). This
completes the proof of Theorem \ref{maintheor}.


\section{Nonmonotonic curves satisfying condition (G)}
\setcounter{equation}{0}

Now consider not necessarily  monotonic curves with fixed Young
diagram $D$ and reduced Young diagram $\Delta$, satisfying condition
(G) (see subsection 3.3). For such curves the canonical complements
$V_i(t)$ to $\Bigl(\Lambda_{(p_{i})}\Bigr)^{(1)}(t)$
 in
$\Lambda_{(p_{i}-1)}(t)$ are defined as well.
 Denote by $\Gamma_{i}^+$ and $\Gamma_{i}^-$ the
positive and the negative index of the quadratic form
$\dot\Lambda(t)|_{(\Lambda_{(p_i-1)})^{(p_i-1)}(t)}$ and let
$r_{i}^+=\Gamma_{i}^+-\Gamma_{i-1}^+$ and
$r_{i}^-=\Gamma_{i}^--\Gamma_{i-1}^-$. Actually the numbers
$r_{i}^+$ and $r_{i}^-$ are equal to
the positive and negative inertia index of the canonical quadratic
forms $Q_{i,t}$ on $V_i(t)$ . These numbers do not depend on $t$ and
they will be called  \emph {the $i$th positive inertia index  and
the $i$th negative inertia index of the curve $\Lambda(t)$} respectively.
Similarly to Definition \ref{normframe} one can define the normal
(quasi-normal) moving frame for a curve in a Lagrange Grassmannian,
satisfying condition (G). The only modification comparing to this
definition is that one should replace the second line in the
structural equation (\ref{structeq}) by $E_a'=I_{r_i^+, r_i^-}
F_a(t)$, $a\in \mathcal F_1\cap \Upsilon_i$, where $r_i^+$ and
$r_i^-$ are the $i$th positive and negative inertia indices of the
curve $\Lambda(t)$,
and the matrix $I_{r_i^+, r_i^-}$ is the diagonal $(r_i^++
r_i^-)\times (r_i^++ r_i^-)-$matrix such that its first $r_i^+$
diagonal entries are equal to $1$ and others are equal to $-1$.
Continuing the normalization procedure by complete analogy with
subsections 3.4-3.6 with obvious modifications, one gets the
following generalization of Theorem \ref{maintheor} to nonmonotonic
curves satisfying condition (G):

\begin{theor}
\label{maintheor1} For any curve $\Lambda(t)$ with the Young diagram
$D$ in the Lagrange Grassmannian, satisfying condition (G), there
exists a normal moving frame $\bigl(\{E_a(t)\}_{a\in
\Delta},\{F_a(t)\}_{a\in \Delta}\bigr)$. A moving frame
$\bigl(\{\tilde e_\alpha(t)\}_{\alpha\in D},\{\tilde
f_\alpha(t)\}_{\alpha\in D}\bigr)$
is a normal moving frame of the curve $\Lambda(\cdot)$ if and only
if for any $1\leq i\leq d$ there exists a constant matrix $U_i\in
O(r_i^+, r_i^-)$ such that for all $t$
\begin{equation}
\label{u1} \widetilde E_a(t)=E_a(t)U_i,\quad \widetilde
F_a(t)=F_a(t)I_{r_i^+, r_i^-}U_i I_{r_i^+, r_i^-}, \quad \forall
a\in \Upsilon_i,
\end{equation}
where $r_i^+$ and  $r_i^-$ are the $i$th positive and the negative
inertia indices of the curve $\Lambda(t)$.
\end{theor}

Further, take a Young diagram $D$, as before,  and fix a tuple  of
nonnegative integers $\{r_i^-\}_{i=1}^d$ such that $0\leq r_i^-\leq
r_i$ for any $1\leq i\leq d$. Let $\mathfrak Q_D$ be the quiver,
defined in subsection 2.3. A representation of the quiver $\mathfrak
Q_D$ will be called \emph{compatible  with the Young diagram $D$ and
the tuple $\{r_i^-\}_{i=1}^d$}, if for any $1\leq i\leq d$ the space
of the representation corresponding to the vertex $\Upsilon_i$ is a
$r_i$-dimensional pseudo-Euclidean space with negative inertia index
$r_i^-$  and the linear mappings $\mathcal R(a,b)$ of the
representation corresponding to arrows $(a,b)$ satisfy the following
relations: $\mathcal R(a,b)^*=\mathcal R(b,a)$ and $\mathcal
R\bigl(a,r(a)\bigr)$ is antisymmetric w.r.t. the corresponding
pseudo-Euclidean structure. Then by complete analogy with Theorem
\ref{size1} we have

\begin{theor}
\label{size2} For the given one-parametric family $\Xi(t)$ of
representations of the quiver $\mathfrak Q_D$ compatible with the
Young diagram $D$ with $|D|$ boxes and the tuple  of nonnegative
integers $\{r_i^-\}_{i=1}^d$ there exists the unique, up to a
symplectic transformation, curve $\Lambda(t)$, satisfying condition
(G), in the Lagrange Grassmannian of $2|D|$-dimensional symplectic
space with the Young diagram $D$ such that the quiver of curvatures
of
$\Lambda(t)$ is isomorphic to $\Xi(t)$ and its $i$th negative
inertia index
is equal to $r_i^-$ for any
$1\leq i\leq d$.
If, in addition,  all rows of
$D$ have different length,
then given a tuple of smooth functions
$\{\rho_{a,b}(t):(a,b)\in\Delta\times \Delta, (a,b) \,\,\text{is an
essential pair} \}$ there exists the unique, up to a symplectic
transformation, curve $\Lambda(t)$, satisfying condition (G), in the
Lagrange Grassmannian of $2|D|$-dimensional symplectic space with
the Young diagram $D$ such that for any essential pair
$(a,b)\in\Delta\times \Delta$ and any $t$ its $(a,b)$-curvature at
$t$ coincides with $\rho_{a,b}(t)$ and its $i$th negative inertia
index
is equal to $r_i^-$ for any $1\leq
i\leq d$.
\end{theor}

\section{Consequences for differential geometry of geometric
structures on manifolds} \setcounter{equation}{0}

Let $\mathcal V$ be a geometric structure on a manifold $M$, as in
the Introduction, and  $H$ be its maximized Hamiltonian, which is
smooth on an open subset
of $T^*M$. Assume that
the point $\lambda \in T^*M$ satisfies: $H(\lambda)> 0$,
$dH(\lambda)\neq 0$, and the germ of the Jacobi curve $J_\lambda(t)$
at $t=0$
has Young diagram $D$ with the reduced diagram $\Delta$ and with
$p_1$ boxes in the first row. Let, as before, $W_\lambda=T_\lambda
\mathcal H_{H(\lambda)}/\{\mathbb R\vec H(\lambda)\}$ be the
symplectic space, where the Jacobi curve $J_\lambda(t)$ lives. The
point $\lambda$ will be called \emph{ $D$-regular} if, in addition
to above,
\begin{equation}
\label{cons0} J_\lambda^{(p_1)}(0)=W_\lambda
\end{equation}
and the germ of the Jacobi curve $J_\lambda(t)$ at $t=0$ satisfies
condition (G). Here for simplicity we will work mainly with
$D$-regular points for some Young diagram $D$. Let
\begin{equation}
\label{cons1} J_\lambda(0)=\oplus_{a\in\Delta} \widetilde{\mathfrak
A}_a(\lambda)
\end{equation}
 be the canonical splitting  of the subspace
$J_\lambda(0)$ (w.r.t. the canonically parameterized curve
$J_\lambda(0)$) and ${\rm proj}_\lambda:T_\lambda{\mathcal
H}_{H(\lambda)}\mapsto W_\lambda$ be the canonical projection on the
factor-space.  Set
\begin{equation}
\label{cons2} {\mathfrak A}_a(\lambda)=({\rm
proj}_\lambda)^{-1}\bigl(\widetilde{\mathfrak
A}_a(\lambda)\bigr)\cap \Pi_\lambda,
\end{equation}
where $\Pi_\lambda$ is the vertical subspace of $T_\lambda \mathcal
H_{H(\lambda)}$, defined by \eqref{pilam}. Taking into account that
${\rm proj}_\lambda$ establishes an isomorphism between
$\Pi_\lambda$ and $J_\lambda(0)$, we get from \eqref{cons1} and
\eqref{cons2}
the following \emph {canonical splitting of  the tangent space
$T_\lambda \bigl(T^*_{\pi(\lambda)}M\bigr)$ to the fiber of $T^*M$
at $\lambda$}:
\begin{equation}
\label{cons3} T_\lambda T^*_{\pi(\lambda)}M= \oplus_{a\in\Delta}
\widetilde{\mathfrak A}_a(\lambda)\oplus {\rm
span}\,\{\epsilon(\lambda)\},
\end{equation}
where $\epsilon$ is the Euler field of $T^*M$, i.e. the
infinitesimal generator of the homotheties of the fibers of $T^*M$.
Besides, each subspace $\mathfrak A_a(\lambda)$ is endowed with the
canonical pseudo-Euclidean structure and the corresponding curvature
mappings between the subspaces of the splitting are intrinsically
related to the geometric structure $\mathcal V$.

 Further, let
\begin{equation}
\label{cons4} {\rm Hor}(\lambda)=({\rm
proj}_\lambda)^{-1}\bigl(J_\lambda^{\rm trans}(0)\bigr),
\end{equation}
where $J_\lambda^{\rm trans}(0)$ is the subspace corresponding to
the canonical complementary curve to the Jacobi curve $J_\gamma$ at
$t=0$. Then ${\rm Hor}(\lambda)$ is transversal to the tangent space
$T_\lambda\bigl(T_{\pi(\lambda)}^*M\bigr)$ to the fiber of $T^*M$ at
$\lambda$. So, if for some diagram $D$ the set $U$ of its regular
$D$-points is open in $T^*M\backslash H^0$, then
 for any $q\in \pi(U)$ the subsets $T^*_qM\cap U$ of the
linear space $T^*_qM$ is endowed with very rich additional
structures: at each point $\lambda\in T^*_qM\cap U$  there is the canonical
splitting of tangent spaces (smoothly depending on $\lambda$) such
that the subspaces of the splitting are parameterized by the
superboxes of the reduced diagram $\Delta$, the dimension of each
subspace is equal to the size of the corresponding superbox, these
subspaces are endowed with the canonical pseudo-Euclidean
structures, and
the canonical linear mappings between these subspaces (i.e. the
$(a,b)$-curvature mappings) are defined. Besides the distribution of
``horizontal'' subspaces $\rm{Hor}(\lambda)$ defines the
\emph{connection on $U\subset T^*M$, canonically associated with
geometric structure $\mathcal V$}.

In the case of sub-Riemannian structures the Hamiltonian $H^2$ is
nonnegative quadratic form on the fibers. First it implies the
monotonicity of the corresponding Jacobi curves. Further assume that
in this case  relation \eqref{cons0} holds for some $\lambda$ and
$p_1$. Then there is a neighborhood $U$ of $\pi(\lambda)$ in $M$ and
an open and dense subset $\mathcal O$ of $U$ that satisfies the
following property: for any $\tilde q\in \mathcal O$ there exists a
neighborhood $\widetilde U\in\mathcal O$ and  a Young diagram $D$
such that for each $\hat q\in \widetilde U$ the intersection of the
set of its $D$-regular points with $T_{\hat q}^*M$ is an nonempty
Zariski open subset of $T_{\hat q}^*M$.
Besides, if one works with the Hamiltonian $H^2$ instead of $H$, the
canonical splitting, the canonical Euclidean structures on the
subspaces of the splitting, the curvature mappings, and the
canonical connection above depend rationally on points of the fibers
of $T^*M$. So, to any sub-Riemannian metric satisfying assumptions
above one can assign very rigid additional structures on $T^*M$.

Condition \eqref{cons0} has the following equivalent description in
terms of the extremal $e^{t\vec H}\lambda$. Projections of the
Pontryagin extremals to the base manifold $M$ are called extremal
trajectories. Conversely, an extremal projected to the given
extremal trajectory is called its {\it lift}. From the Pontryagin
Maximum Principle it follows that the set of all lifts of given
extremal trajectory can be provided with the structure of linear
space. The dimension of this space is called {\it corank of the
extremal trajectory}. It turns out that if condition \eqref{cons0}
holds, then corank  of the extremal trajectory $\pi(e^{t\vec
H}\lambda)$ is equal to $1$. Conversely,  if corank  of the extremal
trajectory $\pi(e^{t\vec H}\lambda)$ is equal to $1$, then
$J_{e^{tH}\lambda}^{(p_1(t))}(0)=W_{e^{tH}\lambda}$ for $t$ from
generic set. Note also that if corank of the extremal trajectory is
greater than 1, then this extremal trajectory is the projection of a
so-called abnormal extremal (a Pontryagin extremal living on zero
level set of the corresponding Hamiltonian).
\medskip

{\bf Conjecture} \emph{(private communication with Andrei Agrachev
and Tohru Morimoto). Any sub-Riemannian metric on a completely
nonholonomic vector distribution has at least one corank 1 extremal
trajectory or, equivalently, not all extremal trajectories of it
are projections of abnormal extremals.}
\medskip

If the conjecture is true, then the construction above can be
implemented for any sub-Riemannian metric on  a completely
nonholonomic vector distribution.

In the case of a Riemannian metric the canonical connection above
coincides with the Levi-Civita connection, the reduced Young diagram
of Jacobi curves consists of only one superbox, and
the corresponding curvature mapping can be identified with the part
of Riemannian curvature tensor  appearing in the classical Jacobi
equation for Jacobi vector fields along the Riemannian
geodesics(\cite{agrgam1}). In particular, the whole Riemannian
curvature tensor can be recovered from it. In general case the
relation between $(a,b)$-curvature mappings and the curvature tensor
of the canonical connection is subject for further study.

Finally, if the Jacobi curve $J_\lambda(t)$
has Young diagram $D$ with $p_1$ boxes in the first row such that
$J_\lambda^{(p_1)}(t)\subsetneq W_\lambda$, then using Remark
\ref{loss}, one can make analogous construction on the space
$J_\lambda^{(p_1)}(0)/(J_\lambda^{(p_1)}(0))^\angle$.
\medskip

{\bf Acknowledgements} We would like to thank professor Andrei
Agrachev for his constant attention for this work and for advising
us to use the language of quivers in formulation of our results on
complete system of invariants. We are also very greatful to
professors Boris Doubrov and Joseph Landsberg for useful
discussions, which helped very much in improving  the exposition.


\end{document}